\documentclass[journal]{IEEEtran}
\usepackage{cite}
\usepackage{amssymb}
\usepackage{multirow}
\usepackage{subfigure}
\usepackage{hhline}
\usepackage{lipsum}
\usepackage{multicol}
\usepackage[flushleft]{threeparttable}
\usepackage{footnote}
\makesavenoteenv{tabular}
\makesavenoteenv{table}
\usepackage[usenames, dvipsnames]{color}
\ifCLASSINFOpdf
   \usepackage[pdftex]{graphicx}
   \DeclareGraphicsExtensions{.pdf,.jpeg,.png}
\else
   \usepackage[dvips]{graphicx}
   \DeclareGraphicsExtensions{.eps}
\fi
\usepackage[cmex10]{amsmath}

\usepackage{algorithm}
\usepackage{algorithmicx}
\usepackage{algpseudocode}
\algnewcommand\algorithmicinput{\textbf{Input:}}
\algnewcommand\Input{\item[\algorithmicinput]}
\algnewcommand\algorithmicoutput{\textbf{Output:}}
\algnewcommand\Output{\item[\algorithmicoutput]}
\algnewcommand\algorithmicstep {\textbf{Step}}
\algnewcommand\Step [1]{\item[\algorithmicstep\ \textbf{#1:}]}
\algnewcommand\algorithmicministep {\emph{(}}
\algnewcommand\MiniStep [1]{\item[\algorithmicministep\emph{#1)}]}
\usepackage{dblfloatfix}
\begin{document}
%
\title{The Normalized Singular Value Decomposition of Non-Symmetric Matrices Using Givens fast Rotations}
%
%
%

\author{Ehsan~Rohani,
        Gwan~S.~Choi, Mi~Lu
\thanks{Ehsan Rohani (ehsanrohani@tamu.edu), Gwan S. Choi (gchoi@ece.tamu.edu), and Mi Lu (mlu@ece.tamu.edu)are with Department of Electrical and Computer Engineering, Texas A\&M University, College Station, Texas 77843.}}

%
%

\markboth{Journal of \LaTeX\ Class Files,~Vol.~13, No.~9, September~2014}%
{Shell \MakeLowercase{\textit{et al.}}: Bare Demo of IEEEtran.cls for Journals}
%



\maketitle

\begin{abstract}
In this paper we introduce the   algorithm and the fixed point hardware to calculate the normalized singular value decomposition of a non-symmetric matrices using Givens  fast (approximate) rotations. This algorithm only uses the basic combinational logic modules such as adders, multiplexers, encoders, Barrel shifters (B-shifters), and comparators and does not use any lookup table. This method in fact combines the iterative properties of singular value decomposition method and CORDIC method in one single iteration.
The introduced architecture is a systolic architecture that uses two different types of processors, diagonal and non-diagonal processors. The diagonal processor calculates, transmits and applies the horizontal and vertical rotations, while the non-diagonal processor uses a fully combinational architecture to receive, and apply the rotations. The diagonal processor uses priority encoders, Barrel shifters, and comparators to calculate the rotation angles. Both processors use a series of adders to apply the rotation angles.
The design presented in this work provides $2.83\sim649$ times better energy per matrix performance compared to the state of the art designs. This performance achieved without the employment of pipelining; a better performance advantage is expected to be achieved employing pipelining.
\end{abstract}
\section{Introduction}
Emerging technologies require a high performance, low power, and efficient solution for singular value decomposition (SVD) of matrices. 
A High performance and throughput hardware implementation of SVD is necessary in applications such as linear receivers for 5G MIMO telecommunication systems \cite{SlowestDscentMethod}, various real-time applications \cite{ahmedsaid2004accelerating}, classification in genomic signal processing \cite{genomic}, and learning algorithm in active deep learning \cite{DL}.
Matrix decomposition, and more specifically, SVD is also the most commonly used DSP algorithm and often is the bottle neck of various computationally intensive algorithms. For instance this is paramount for some the resent studies including performance analysis of adaptive MIMO transmission in a cellular system \cite{ref1}, image compression \cite{ref2}, and new image processing techniques of face recognition \cite{ref4}.\par

The design of SVD arithmetic unit has been vastly investigated by the researchers. For effective implementation of SVD in the hardware, \cite{blv1983} presents BLV algorithm with a systolic architecture. This architecture uses a set of diagonal processors (DP) and a set of non-diagonal processors (NDP). The DP processor calculates, applies and transmits the horizontal ($\theta_H$) and vertical ($\theta_V$) rotation angles while NDP applies the received rotation angles.
To calculate the division, square root, and multiplication required for this method Cavallaro uses the coordinate rotation digital computer algorithm (CORDIC) \cite{cavallaro1988cordic}.
CORDIC algorithm is mainly credited to Jack Volder. This algorithm was originally introduced for solving the problem of real-time navigation \cite{andraka1998survey}. CORDIC algorithm has been used in different math coprocessor \cite{DCT}, digital signal processors \cite{Processor}, and software defined radios \cite{SDR}.  Different implementation of CORDIC algorithm have been introduced, including but not limited to Higher Radix CORDIC algorithms \cite{arch1}, Angle Recoding methods \cite{arch2}, Hybrid and Coarse-Fine Rotation CORDIC \cite{arch3}, Redundant-Number-Based CORDIC implementation \cite{arch4}, Pipelined CORDIC architecture \cite{arch5}, and Differential CORDIC algorithm \cite{arch6}. A relatively comprehensive review of CORDIC algorithm is presented in \cite{50years}. To reduce the implementation complexity of SVD, different optimization have been proposed. Delsome proposed double rotations to avoid square roots and divisions for scaling \cite{delosme1992bit}.\par 
In 1991 Gotze introduces an algorithm that combines the inner-iterations (CORDIC iterations) with outer-iterations of BLV algorithm (sweep). This method instead of calculating the accurate rotations, calculates the fast rotation (Givens fast rotations also known as Givens approximate rotations) angles which is equivalent to one iteration of CORDIC algorithm \cite{gotze1991parallel}, \cite{gotze1993efficient}. As the result this hardware does not require any look-up table, and calculating the rotation angles requires only one clock cycle. However this method has the following disadvantages:\par
\begin{enumerate}
\item The hardware implementation requires floating-point arithmetic.\par
\item  The algorithm only works for symmetric matrices.\par
\item Using this method the BLV algorithm loses its quadratic convergence speed (quadratic to number of sweeps).\par
\item The proposed algorithm does not provide the "Normalized" results.\par
\end{enumerate}
In this work we address the disadvantages of Givens fast rotations. The proposed hardware does not require the floating-point arithmetic, and as the result, it does not need the pre-processing (alignment of exponents) and post-processing (renormalization of matrices) blocks mentioned in implementations of Givens fast rotation. The algorithm that we proposed is able to handle symmetric as well as non-symmetric matrices (calculating the rotations for non-symmetric matrices require the calculation of two intermediate variable angles, proximate calculation of these two intermediate variables makes the calculation of rotation angles challenging and we were able to offer an adaptive solution for it). Also, the proposed algorithm provides the "Normalized" results. 
The rest of this work is organized as follows: First we introduce the method to merge the NSVD algorithm \cite{blv1983} with Givens fast rotations and Delsome double rotations method \cite{delosme1989cordic}. The result (ERNSVD algorithm) is similar to the Gotze work in \cite{gotze1993efficient} to find the Eigenvalues of a symetric matrix exept the algorithm is able to calculate the decomposition of a non-symetric matrix. In addition we introduce a method that directly calculates the horizontal and vertical fast rotations using the Forsythe and Henrici SVD (FHSVD) algorithm and called it expedite rotations SVD (ERFHSVD) \cite{blv1983}. We present a possible hardware implementation and provisions that make the fixed point implementation of these algorithm possible in section \ref{HWI}. The hardware implementation is a design for decomposing a $2 \times 2$ matrix, as the basic building block for decomposition of matrices with larger size. In section \ref{Complexity} the complexity of the suggested design is estimated based on two factors, resource requirement and the critical path delay.
 Finally we conclude this Work in section \ref{Concl}.\par
\section{Algorithm of calculating fast rotation angles for non-symmetric matrices}
\label{Algorithm}
This work is inspired by the Normalized SVD (NSVD) and FHSVD algorithm presented in \cite{blv1983}.
We use the Givens fast rotations presented in \cite{gotze1993efficient} and double rotation by Delsome \cite{delosme1992bit} to reduce the implementation complexity. 
In this section we briefly review the bases of our inspiration (NSVD, FHSVD algorithms, double rotations, and fast rotations) and introduce our measures (approximation error and norm of off-diagonal elements) for evaluation and comparison of different methods through this work. This section holds four subsections containing our three proposed algorithms as well as their comparison and also study of the relaxing the boundary conditions which will result in reduction in hardware complexity.\par
The fast rotation algorithm tries to find the closest angle ($\hat{x}$) to any rotation angle ($x$) such that: first, $|\hat{x}| \leq \frac{\pi}{4}\quad$ and also $|tan(\hat{x})|=2^{-l}; \quad l \geq 0$. The fast rotations do not symmetrize or diagonalize the matrix in one rotation, in fact they generate a more symmetric matrix with each rotation or they reduce the off-diagonal norm of the matrix; using this method the quadratic convergence properties of SVD algorithm is lost. 
We use the Delsome proposed double rotations method which for any rotation angle uses $\bar{\gamma}  = \frac{\gamma}{2}$ as rotation angles and applies the rotations twice \cite{delosme1989cordic}. This technique will eliminate the need for the calculation of square root function and division for scaling. 
This means adding one to $l$ and applying two rotations with the angle equal to $arctan(2^{-l+1})$. We use $\tilde{\ }$ for noting any angle or its tangent when Delsome double rotation and Gotze approximation are applied. \par

We assume the $2 \times 2$ matrix $A$ and the decomposed matrix is defined as in (\ref{eq:1}) for the review of NSVD and FHSVD algorithms.  

\begin{equation}
\mathbf{A}=
\begin{bmatrix}
a & b\\ 
c & d
\end{bmatrix}
= \mathbf{U} \times \pmb{\Sigma} \times \mathbf{V}^T 
\label{eq:1}
\end{equation}

The NSVD algorithm calculates the first rotation to symmetrize the matrix using (\ref{eq:sym}, and then using (\ref{eq:symMult}) generates the symmetric matrix. NSVD then uses (\ref{eq:diag}) to find the diagonalizing rotations.
\begin{equation}
 \rho = tan^{-1}( \frac{c+b}{d-a} )
\label{eq:sym}
\end{equation}
\begin{equation}
 \mathbf{B}=\mathbf{R}_{\rho} \times \mathbf{A}=\begin{bmatrix}cos(\rho) & sin(\rho) \\-sin(\rho) & cos(\rho) \end{bmatrix}\begin{bmatrix}a & b \\c & d \end{bmatrix}    = \begin{bmatrix}p & q \\q & r \end{bmatrix}
\label{eq:symMult}
\end{equation}
\begin{equation}
 \phi = tan^{-1}( \frac{2q}{q-p} )
\label{eq:diag}
\end{equation}
\begin{equation}
\begin{aligned}
  \pmb{\Sigma}=\mathbf{R}_{\phi}^T \times \mathbf{B} \times \mathbf{R}_{\phi} = \begin{bmatrix}d_1 & 0 \\0 & d_2 \end{bmatrix} \\ =\begin{bmatrix}cos(\phi) & sin(\phi) \\-sin(\phi) & cos(\phi) \end{bmatrix}^T\begin{bmatrix}p & q \\q & r \end{bmatrix}  \begin{bmatrix}cos(\phi) & sin(\phi) \\-sin(\phi) & cos(\phi) \end{bmatrix} 
 \end{aligned}
\label{eq:diagMult}
\end{equation}

The FHSVD algorithm first calculates the $\alpha$ and $\beta$ using (\ref{eq:2}) and (\ref{eq:3}) as the intermediate values, and using (\ref{eq:4}) and (\ref{eq:5}) the algorithm calculates the horizontal ($\Theta$) and vertical ($\theta$) rotations. In the last step FHSVD algorithm swaps the values of $sin$ and $cos$ and change the sign of these values if needed to make sure $d_1 \geq d_2$.

\begin{equation}
\label{eq:2}
 \alpha = tan^{-1}( \frac{c+b}{d-a} )
\end{equation}

\begin{equation}
\label{eq:3}
 \beta = tan^{-1}( \frac{c-b}{d+a} )
\end{equation}

\begin{equation}
\label{eq:4}
 \Theta = \frac{\alpha+\beta}{2}
\end{equation}

\begin{equation}
\label{eq:5}
 \theta = \frac{\alpha-\beta}{2}
\end{equation}

\begin{equation}
\begin{aligned}
 \boldsymbol {\Sigma} = \mathbf{R}_{\theta}^T \times \mathbf{A} \times \mathbf{R}_{\Theta} = \begin{bmatrix}d_1 & 0 \\0 & d_2 \end{bmatrix} \\ =\begin{bmatrix}cos(\theta) & sin(\theta) \\-sin(\theta) & cos(\theta) \end{bmatrix}^T\begin{bmatrix}a & b \\c & d \end{bmatrix}  \begin{bmatrix}cos(\Theta) & sin(\Theta) \\-sin(\Theta) & cos(\Theta) \end{bmatrix} 
 \end{aligned}
\label{eq:diagMult}
\end{equation}

\cite{gotze1993efficient} offers a floating point implementation and uses the exponent bits to calculate $l$. In fix point representation we introduce (\ref{eq:7}) to replace the exponent bits value of the original algorithm.
This will reduce the function $exp_2(x)$ to a priority encoder applied on $x$, assuming $x$ is an integer number. We will discuss this in more details in \ref{HWI}. Equation system (\ref{eq:6}) is equal to checking the most significant bit (MSB) in two's-complement representation.

\begin{equation}
\label{eq:6}
Sign(x)= \begin{cases}1 &\text{ if } x  \geq 0\\-1 &\text{ if } x  <  0\end{cases}
\end{equation}  

\begin{equation}
\label{eq:7}
(exp_2(x),v)= \begin{cases}
(0,0) & \text{ if } x=0 \\ 
(\lfloor {log_{2}}^{\left | x \right |} \rfloor ,1) & \text{ if } x\neq 0
\end{cases}
\end{equation}

There are two measures that are normally used in different Jacobi based decomposition. 

\begin{enumerate}
\item The approximation error $\mid d \mid$ in \cite{gotze1993efficient} for a $2 \times 2 $ symmetric matrix defined as the absolute value of off diagonal element before and after the application of $k^{th}$ rotation (\ref{eq:absd}). The  smaller $\mid d \mid$ shows a more accurate approximation.

\begin{equation}
\label{eq:absd}
\mid d \mid= \mid \frac{q_{(k+1)}}{q_k} \mid
\end{equation} 

This is applied to the $2 \times 2$ matrix, and $\mid d \mid _{Max}<1$ is called error bound and is one of the two conditions to assure the convergence of the algorithm. The other condition is to keep the orthogonality of the rotation matrix. We plan to apply this measures to non-symmetric matrices so we extend the equation (\ref{eq:absd}) to $\mid\mid D \mid\mid$. Note that $\mid d \mid$ is the same as $\mid\mid D \mid\mid$ if the matrix is symmetric.

\begin{equation}
\label{eq:absD}
\mid\mid D \mid\mid=  \frac{\sqrt{b^2_{k+1}+c^2_{k+1}}}{\sqrt{b^2_{k}+c^2_{k}}}
\end{equation} 
The approximation error (as the dependent variable) is normally measured for different values of $\tau$ (as the independent variable) as defined in ($\ref{eq:tau}$). 
\begin{equation}
\label{eq:tau}
\tau=  \frac{r-p}{2q}
\end{equation} 

To extend equation (\ref{eq:tau}) to an independent variable applicable to non-symmetric matrices, we define $\tau_1$ and $\tau_2$ in (\ref{eq:tau1}) and (\ref{eq:tau2}) accordingly. $\frac{1}{\tau_2}=0$ and $\tau_1=\tau$ if the matrix is symmetric.
     
\begin{equation}
\label{eq:tau1}
\tau_1=  \frac{d-a}{b+c}
\end{equation} 

\begin{equation}
\label{eq:tau2}
\tau_2=  \frac{d+a}{b-c}
\end{equation} 

\item The Norm of the off-diagonal elements of a matrix versus (Vs.) the number of sweeps is the second metric used for measuring the quality of any fast rotation methods as well as the comparison of the diagonalization speed in different methods. In \cite{gotze1993efficient} this value is calculated for a $70 \times 70$ matrix. The elements of matrix are randomly generated numbers of normal distribution. We run the same test for 100 times and calculate the RMS (root mean square) of the off-diagonal norms ($RMS_{ODN}$) to keep the comparability and also keep the results accurate.
\end{enumerate} 

\begin{figure}
\centering
\includegraphics[width=0.49\textwidth,keepaspectratio]{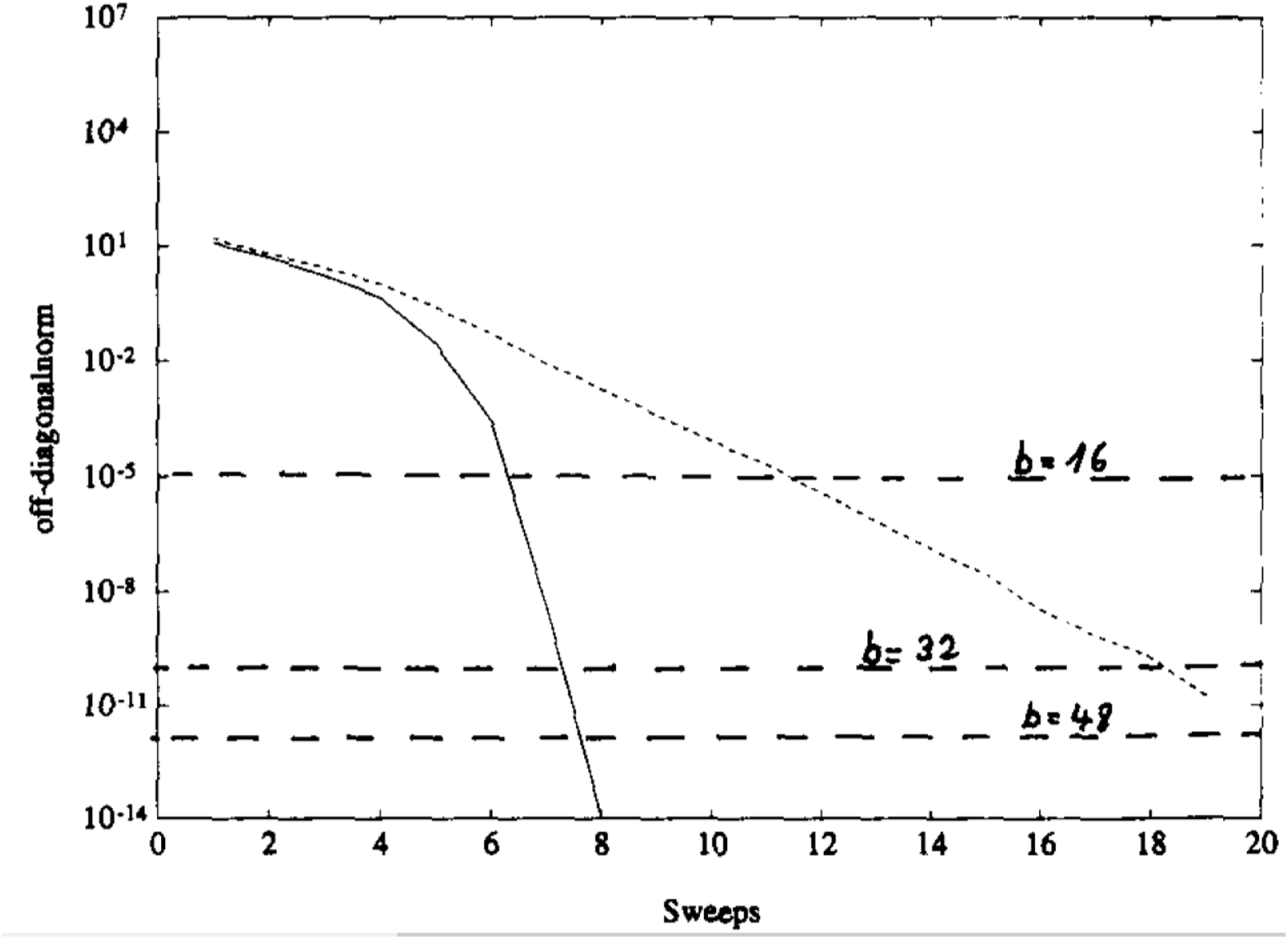}{}
\caption{Off-diagonal Norm vs. number of sweeps \cite{gotze1993efficient}}
\label{fig:ODN}
\end{figure}

Figure~\ref{fig:ODN} shows the Off-diagonal Norm of a $70 \times 70$ vs. number of sweeps for the fast rotations (dotted curve) and the original Givens rotations (solid line). The dashed lines show the accuracy achievable by different number of bits. This figure provided in \cite{gotze1993efficient} is to support the fast rotations method with the following explanation: While achieving the accuracy of 16 bits or better; the original method requires seven iterations versus twelve iterations in the fast rotation method. One must realize that the calculation of 16-bit exact sine and cosine values requires more complexity. As an example if the sine and cosine are implemented using CORDIC method, it requires 16 inner iteration to calculate these values. In the next two subsection we will present two algorithms that are basic blocks of eminent rotations NSVD (ERNSVD).

\subsection{Symmetrizing Algorithm}
\label{SA}
Algorithm~\ref{symalg} is the approximation to the rotation angle achieved from equation (\ref{eq:sym}). Boundaries in Step~2 of the algorithm are based on the suggestion in \cite{gotze1991parallel} that guaranties $|d| \leq \frac{1}{3}$. Authors in \cite{gotze1993efficient} explain that using Delsome double rotations this limit does not hold any more (this method will guaranties $|d| \leq \frac{7}{12}$); however, the authors explain that this is not an issue since the approximation error merges to the original bounds when the rotation angles get smaller. It should be mentioned that the bounds to approximation error is still less than one.

\begin{algorithm}
\caption{Algorithm to Calculate the Symmetrizing Fast Rotations for Non-Symmetric Matrices.}
\label{symalg}
\begin{algorithmic}
\Input {$\mathbf{A}$.}
\Output	{Rotation Matrix $\mathbf{R}_{\tilde{\rho}}$.}
\Step{1} Calculate the initial values:
$$S_{N}=Sign(b-c)$$
$$S_{D}=Sign(d+a)$$
$$N=|b-c|$$
$$D=|d+a|$$
$$K=exp_2(D)-exp_2(N)$$
\Step{2} Calculate $l_{\tilde{\rho}}$ using following case statement:
$$l_{\hat{\rho}}=\begin{cases}
K+1 & \text{ if } 1.5\times D> (2^{K+1}-2^{-K})\times N \\ 
K-1 & \text{ if } 1.5\times D< (2^{K}-2^{-K+1})\times N \\ 
K & \text{ default } 
\end{cases}$$
$$l_{\tilde{\rho}}=max(l_{\hat{\rho}}+1,1)$$
\Step{3} Calculate the $(c,s)$ pair using the following case statement:
$$\tilde{t}= 2^{-l_{\tilde{\rho}}}$$
$$(c,s)=\begin{cases}
(1,0) & \text{ if } N=0\\ 
(0,1) & \text{ if } D=0\\ 
\frac{1}{1+\tilde{t}^2} \times (1-\tilde{t}^2,2 \times S_{N} \times S_{D} \times \tilde{t}) & \text{ default } 
\end{cases}$$
\Step{4} Calculate the Symmetrizing $\mathbf{R}_{\tilde{\rho}}$:
$$\mathbf{R}_{\tilde{\rho}}=
\begin{bmatrix}
c & s\\ 
-s & c
\end{bmatrix}
$$
\end{algorithmic}
\end{algorithm}

\subsection{Diagonalizing Algorithm}
\label{DA}

\begin{figure}
\centering
\includegraphics[width=0.49\textwidth,keepaspectratio]{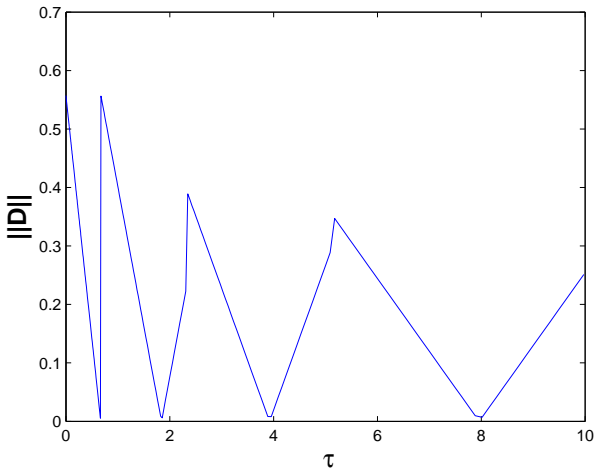}{}
\caption{$\mid\mid D \mid\mid$ vs. $\tau$ for symmetric matrix when Algorithm~\ref{diagalg} is applied.}
\label{fig:DVsT}
\end{figure}

\begin{figure}
\centering
\includegraphics[width=0.49\textwidth,keepaspectratio]{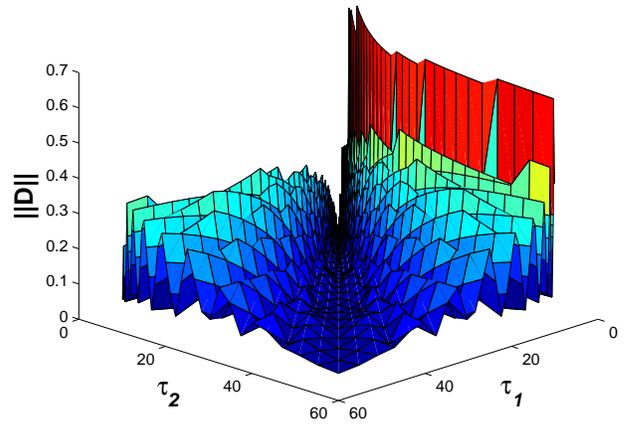}{}

\caption{$\mid\mid D \mid\mid$ vs. $\tau_1$ and $\tau_2$ when Algorithm~\ref{diagalg} is applied for non-symmetric matrix with ideal symmetrization.}
\label{fig:DVsT1T2Idl}
\end{figure}

Algorithm~\ref{diagalg} is the approximation to the rotation angle achieved from equation (\ref{eq:diag}).
 Figure~\ref{fig:DVsT} demonstrates the $\mid\mid D \mid\mid$ vs. $\tau$ when Algorithm~\ref{diagalg} is applied on a symmetric matrix. Figure~\ref{fig:DVsT1T2Idl} demonstrates the $\mid\mid D \mid\mid$ vs. $\tau_1$ and $\tau_2$ when Algorithm~\ref{diagalg} is applied on a non-symmetric matrix with ideal symmetrizing rotations. The fast rotations NSVD (FRNSVD) algorithm is based on applying the fast symmetrizing and fast diagonalizing algorithms on a given matrix. \par

\begin{figure}
\centering
\includegraphics[width=0.49\textwidth,keepaspectratio]{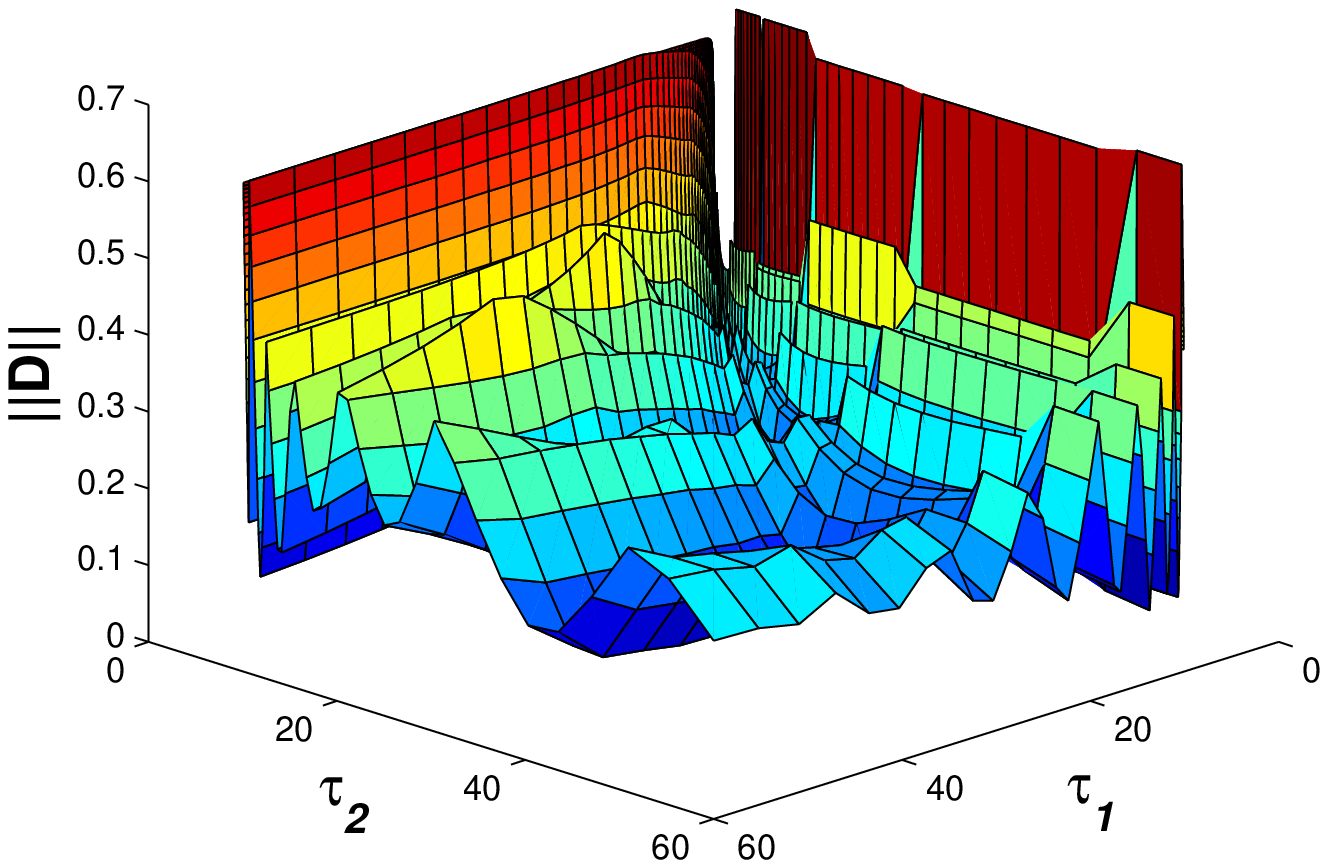}{}
\caption{$\mid\mid D \mid\mid$ vs. $\tau_1$ and $\tau_2$ for non-symmetric matrix for FRNSVD.}
\label{fig:DVsT1T2}
\end{figure}

Figure~\ref{fig:DVsT1T2} demonstrates the $\mid\mid D \mid\mid$ vs. $\tau_1$ and $\tau_2$ when Algorithm~\ref{symalg} and Algorithm~\ref{diagalg} is applied on a non-symmetric matrix. The approximation errors are typically higher in this method compared to Figure~\ref{fig:DVsT1T2Idl}. This increase in the approximation error is expected since the approximation error of two algorithms can boost the total approximation error.
Figure~\ref{fig:NormCompRelax} compares the $RMS_{ODN}$ for NSVD, FRNSVD, and ERNSVD. The error floor that happens in iterations 19th and after in ERNSVD, and FRNSVD are due to the fixed point (32 bit) implementation of the algorithm. 
\begin{algorithm}
\caption{Algorithm to Calculate the fast Diagonalizing Rotations for Symmetric Matrices.}
\label{diagalg}
\begin{algorithmic}
\Input {$\mathbf{B}$.}
\Output	{Rotation Matrix $\mathbf{R}_{\tilde{\phi}}$.}
\Step{1} Calculate the initial values:
$$S_{N}=Sign(b+c)$$
$$S_{D}=Sign(d-a)$$
$$N=|c+b|$$
$$D=|d-a|$$
$$K=exp_2(D)-exp_2(N)$$
\Step{2} Calculate $l_{\tilde{\phi}}$ using following case statement:
$$l_{\hat{\phi}}=\begin{cases}
K+1 & \text{ if } 1.5\times D> (2^{K+1}-2^{-K})\times N \\ 
K-1 & \text{ if } 1.5\times D< (2^{K}-2^{-K+1})\times N \\ 
K & \text{ default } 
\end{cases}$$
$$l_{\tilde{\phi}}=max(l_{\phi}+2,1)$$
\Step{3} Calculate the $(c,s)$ pair using the following case statement:
$$\tilde{t}= 2^{-l_{\tilde{\phi}}}$$
$$(c,s)=\begin{cases}
(1,0) & \text{ if } N=0\\ 
(0,S_{N}) & \text{ if } D=0\\ 
\frac{1}{1+\tilde{t}^2} \times (1-\tilde{t}^2,2 \times S_{N} \times S_{D} \times \tilde{t}) & \text{ default } 
\end{cases}$$
\Step{4} Calculate the diagonalizing $\mathbf{R}_{\tilde{\phi}}$ using the following case statement:
$$ \mathbf{R}_{\tilde{\phi}} =\begin{cases}
\begin{bmatrix} c & s \\-s & c \end{bmatrix} & \text{ if } S_N<0\\ \\
\begin{bmatrix} s & c \\c & -s \end{bmatrix} & \text{ else } 
\end{cases}$$

\end{algorithmic}
\end{algorithm}

\subsection{Relaxing the Boundary Conditions}
\label{RBC}
Implementing the case statement in \textbf{Step 2} of both algorithms beside the two comparators requires two Barrel shifters,two adders to calculate the coefficients of \emph{N} and an adder to calculate $1.5 \times D$. We can reduce the complexity of conditions by using (\ref{eq:relax1}). This reduces the complexity of case statement to two comparators, an adder, and a Barrel shifter. \par

We study the effect of relaxing the conditions using both of the measures.  Equation (\ref{eq:relax1}) provides an approximation that is accurate for higher \emph{l}s while it is inaccurate for small \emph{l}s. The solution is in limiting the rotation angles to smaller rotations (changing $l=max(l+1,1)$ to $l=max(l+1,2)$); in another word, the rotation angles are over estimated when the original rotation is closer to $\dfrac{\pi}{4}$. The empirical results of the simulation shows that changing both $l_\rho$ and $l_\phi$ will result in reduced convergence speed. We found that the best combination is achieved by limiting the $l_\rho$ to minimum of two while $l_\phi$ can get any positive integer value. \par 
\begin{equation}
\label{eq:relax1}
l=\begin{cases}
K+1 & \text{ if } 1.5\times D> 2^{K+1}\times N \\ 
K-1 & \text{ if } 1.5\times D< 2^{K}\times N \\ 
K & \text{ default } 
\end{cases}
\end{equation}

\begin{figure}
\centering
\includegraphics[width=0.49\textwidth,keepaspectratio]{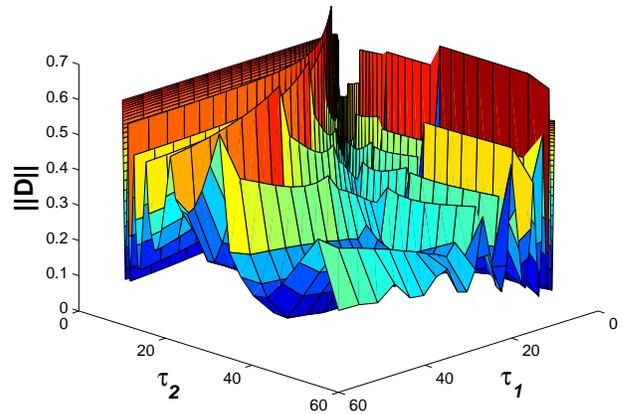}{}
\caption{$\mid\mid D \mid\mid$ vs. $\tau_1$ and $\tau_2$ for non-symmetric matrix for ERNSVD.}
\label{fig:DVsT1T2Simp}
\end{figure}

\begin{figure}
\centering
\includegraphics[width=0.49\textwidth,keepaspectratio]{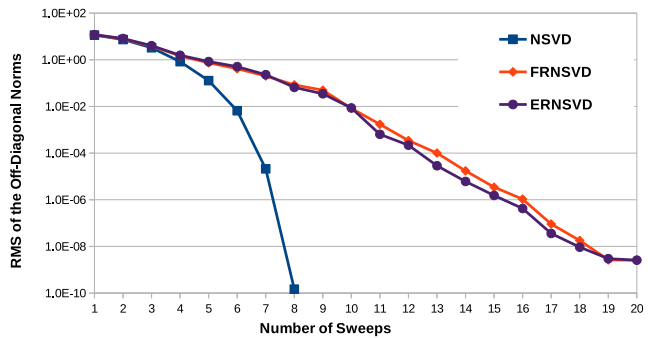}{}
\caption{Comparison of the $RMS_{ODN}$ for NSVD, FRNSVD, ERNSVD algorithms.}
\label{fig:NormCompRelax}
\end{figure}

We called the relaxed version of fast rotations the expedite rotations NSVD (ERNSVD). Figure~\ref{fig:NormCompRelax} compares the $RMS_{ODN}$ for NSVD, FRNSVD and ERNSVD. The $RMS_{ODN}$ are very close for both methods and with ERNSVD being less complex the new boundaries are studied after this point. Figure~\ref{fig:DVsT1T2Simp} demonstrates the $\mid\mid D \mid\mid$ vs. $\tau_1$ and $\tau_2$ for ERNSVD. The FRNSVD has lower approximation error compare to ERNSVD when $\tau_1$ or $\tau_2$ are closer to zero. This is the effect of the empirical change we applied to the boundaries in equation (\ref{eq:relax1}).

\begin{figure*}
\centering
\includegraphics[width=0.98\textwidth,keepaspectratio]{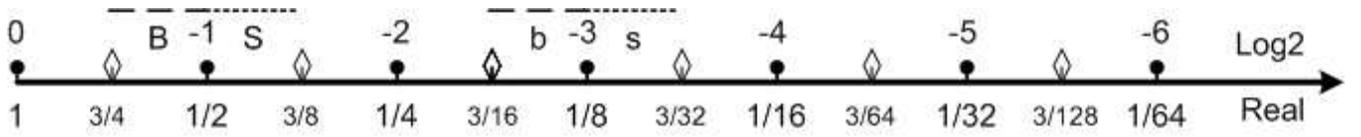}{}

\caption{Axis showing the boundary conditions of each fast rotation angle for relaxed version}

\label{fig:Ax}
\end{figure*}

\subsection{Direct Estimate Algorithm}
\label{DEA}
One major contribution of this work is the algorithm that enable us to directly calculate the fast rotations for non-symmetric matrices based on FHSVD algorithm. Assuming a $2 \times 2$ matrix is defined as in equation (\ref{eq:1}), FHSVD algorithm uses (\ref{eq:2}) and (\ref{eq:3}) to calculate $\alpha$ and  $\beta$; then using those two values and equation (\ref{eq:4}) and (\ref{eq:5}), it calculates the vertical and horizontal Givens rotations. Note that for traditional  rotation algorithm (proposed by Gotze \cite{gotze1993efficient}), the value of $\beta$ is zero,  and we only need to calculate the value of $\alpha$ (Gotze algorithm is for symmetric matrices only).\par 

While Gotze method \cite{gotze1993efficient} is using the same flow as in the calculation of accurate rotation angles and then divides the results by two to take benefit of double rotations, we use $tan(x) \approx x$ and the fact that this approximation is more accurate for smaller angels ($lim_{x \longrightarrow 0}^{tan(x)} = x$), to provide more accurate estimation of rotation angels.\par 

Figure \ref{fig:Ax} shows the boundaries for the relaxed version with a rhombus indicator on the axis while the approximated angles are shown with circle indicators. This shows that as an example if the range of the tangent is between ($\frac{3}{8}$,$\frac{3}{4}$) its approximation will be $\frac{1}{2}$. The notion \textbf{B} or \textbf{b} are used when the original rotation angle is larger than approximations while \textbf{B} is used for the bigger angle between $\alpha$ and $\beta$. The notion \textbf{S} or \textbf{s} are used when the rotation angle is smaller than its approximation while \textbf{L} is used for the bigger angle between $\alpha$ and $\beta$. Knowing if the fast rotation angle is overestimating or underestimating the original rotation angle might requires an extra comparator (depends on the implementation) in angle calculation circuit but it could provide considerable benefit as we will demonstrate in following sections\par

TABLE~\ref{tab:Direct} gives a better understanding of how the rotation angles ($\tilde{\theta}$ and $\tilde{\Theta}$) are decided only based on having an approximation of $\alpha$ and $\beta$. $\tilde{\alpha}$ and $\tilde{\beta}$ are approximations to $\alpha$ and $\beta$ that Delsome double rotation and angle approximation are both applied (In TABLE~\ref{tab:Direct} $\alpha$ and $\beta$ are half the rotation angles achieved from (\ref{eq:2}) and (\ref{eq:3}) to take advantage of tangent properties: $lim_{x \longrightarrow 0}^{tan(x)} = x$). Column $\theta$ in the table shows the possible range of $\theta$ based on the value of $\tilde{\alpha}$, $\tilde{\beta}$, and if the approximation angels are larger or smaller than the original rotation angles. $\tilde{\theta}_{big}$ is the approximation if we use the biggest possible rotation angle in the range, and $\tilde{\theta}_{small}$ is the approximation if we use the smallest possible rotation angle in the range (same notation is applied to $\tilde{\Theta}_{big}$ and $\tilde{\Theta}_{small}$). Choosing the bigger rotation angles in general might result in increased floor level of $RMS_{ODN}$ while it results in a faster convergence rate. Columns $\tilde{\theta}$ and $\tilde{\Theta}$ in the table show the approximation we used in algorithm~\ref{alg1} as an example; however, this does not mean that any application of the algorithm has to use the same numbers. In fact we urge a search on the possibilities based on the application.\par

\begin{table*}[]
\small
\centering
\caption{Different possibilities of choosing the fast rotation angles}
\label{tab:Direct}
\begin{tabular}{|c|c||c|c|c|c||c|c|c|c|}
\hline
 & & & & & & & & &\\
$\tilde{\alpha} \quad and \quad \tilde{\beta}$                                                                                           & Range & $\theta=\alpha-\beta$         & $\tilde{\theta}_{big}$ & $\tilde{\theta}_{small}$ & $\tilde{\theta}$ & $\Theta=\alpha+\beta$ & $\tilde{\Theta}_{big}$ & $\tilde{\Theta}_{small}$ & $\tilde{\Theta}$ \\ \hline
\multirow{4}{*}{\begin{tabular}[c]{@{}c@{}}$\tilde{\alpha}=\frac{1}{2}$\\ $\tilde{\beta}=\frac{1}{2}$\end{tabular}}  & Bb       & $(\frac{-1}{4},\frac{1}{4})$                                &                            $\frac{1}{4} \quad or \quad \frac{-1}{4}$ &                0               &        0               &    $(1,\frac{6}{4})$                     &            1                 &            1                   &                      1 \\ \cline{2-10} 
                                                                                                                         & Bs       & $(0,\frac{3}{8})$     &     $\frac{1}{4}$                      &  0   &                      0   & $(\frac{7}{8},\frac{5}{4})$                        &                  1           &                              1 &         1              \\ \cline{2-10} 
                                                                                                                         & Sb       & $(\frac{-3}{8},0)$                                 &               $\frac{-1}{4}$              &       0                        &        0               & $(\frac{7}{8},\frac{5}{4})$                        &            1                 &           1                    &           1            \\ \cline{2-10}
                                                                                                                         & Ss       &  $(\frac{-1}{8},\frac{1}{8})$                               &            $\frac{1}{8} \quad or \quad \frac{-1}{8}$                 &           0                    &            0           & $(\frac{3}{4},1)$                        &                1            &              1                 &                      1 \\ \hline \hline
\multirow{4}{*}{\begin{tabular}[c]{@{}c@{}}$\tilde{\alpha}=\frac{1}{2}$\\ $\tilde{\beta}=\frac{1}{4}$\end{tabular}}  & Bb       & $(\frac{1}{8},\frac{1}{2})$ & $\frac{1}{2}$                            &    $\frac{1}{8}$                           & $\frac{1}{4}$                      & $(\frac{3}{4},\frac{9}{8})$                        & $1$                            &    $1$                           & $1$                      \\ \cline{2-10} 
                                                                                                                         & Bs       & $(\frac{1}{4},\frac{9}{16})$                                & $\frac{1}{2}$                            & $\frac{1}{4}$                              & $\frac{1}{4}$                      & $(\frac{11}{16},1)$                         & $1$                             &    $\frac{1}{2}$                           & $\frac{1}{2}$                       \\ \cline{2-10} 
                                                                                                                         & Sb       & $(0,\frac{1}{4})$                                &             $\frac{1}{4}$                 &                              0 &          $\frac{1}{4}$            & $(\frac{5}{8},\frac{7}{8})$                        &    $1$                         & $\frac{1}{2}$                              & $\frac{1}{2}$                      \\ \cline{2-10} 
                                                                                                                         & Ss       & $(\frac{1}{8},\frac{5}{16})$                                & $\frac{1}{4}$                            & $\frac{1}{8}$                              & $\frac{1}{4}$                      & $(\frac{9}{16},\frac{3}{4})$                         & $\frac{1}{2}$                            & $\frac{1}{2}$                              & $\frac{1}{2}$                      \\ \hline \hline
\multirow{4}{*}{\begin{tabular}[c]{@{}c@{}}$\tilde{\alpha}=\frac{1}{2}$\\ $\tilde{\beta}=\frac{1}{8}$\end{tabular}}  & Bb       & $(\frac{5}{16},\frac{5}{8})$                                & $\frac{1}{2}$                            & $\frac{1}{4}$                               & $\frac{1}{2}$                      &    $(\frac{5}{8},\frac{15}{16})$                     & $1$                             & $\frac{1}{2}$                               & $\frac{1}{2} $                       \\ \cline{2-10} 
                                                                                                                         & Bs       &  $(\frac{3}{8},\frac{21}{32})$                                & $\frac{1}{2}$                             &    $\frac{1}{2}$                           & $\frac{1}{2}$                       & $(\frac{19}{32},\frac{7}{8})$                         & $1$                            &    $\frac{1}{2}$                           & $\frac{1}{2}$                      \\ \cline{2-10} 
                                                                                                                         & Sb       & $(\frac{3}{16},\frac{3}{8})$                                & $\frac{1}{4}$                            & $\frac{1}{4}$                              & $\frac{1}{2}$                      & $(\frac{1}{2},\frac{11}{16})$                        & $\frac{1}{2}$                            & $\frac{1}{2}$                              & $\frac{1}{2}$                      \\ \cline{2-10} 
                                                                                                                         & Ss       & $(\frac{1}{4},\frac{13}{32})$                                &    $\frac{1}{2}$                         & $\frac{1}{4}$                              &    $\frac{1}{2}$                   & $(\frac{17}{32},\frac{5}{8})$                        & $\frac{1}{2}$                             &    $\frac{1}{2}$                           & $\frac{1}{2}$                       \\ \hline \hline
\multirow{4}{*}{\begin{tabular}[c]{@{}c@{}}$\tilde{\alpha}=\frac{1}{2}$\\ $\tilde{\beta}=\frac{1}{16}$\end{tabular}} & Bb       & $(\frac{13}{32},\frac{11}{16})$                                & $\frac{1}{2}$                             & $\frac{1}{2}$                               & $\frac{1}{2}$                       & $(\frac{9}{16},\frac{27}{32})$                        & $1$                            & $\frac{1}{2}                             $ & $\frac{1}{2}$                      \\ \cline{2-10} 
                                                                                                                         & Bs       & $(\frac{7}{16},\frac{45}{64})$                                & $\frac{1}{2}$                            &    $\frac{1}{2}$                           & $\frac{1}{2}$                       & $(\frac{35}{64},\frac{13}{16})$                        & $1$                             & $\frac{1}{2}$                              & $\frac{1}{2}$                      \\ \cline{2-10} 
                                                                                                                         & Sb       & $(\frac{9}{32},\frac{7}{16})$                                &    $\frac{1}{2}$                         & $\frac{1}{4}$                              & $\frac{1}{2}$                       & $(\frac{7}{16},\frac{19}{32})$                        & $\frac{1}{2}$                            & $\frac{1}{2}$                               & $\frac{1}{2}$                      \\ \cline{2-10} 
                                                                                                                         & Ss       & $(\frac{5}{16},\frac{29}{64})$                                & $\frac{1}{2}$                             & $\frac{1}{4}$                              & $\frac{1}{2}$                      & $(\frac{27}{64},\frac{9}{16})$                        & $\frac{1}{2}$                            & $\frac{1}{2}$                              & $\frac{1}{2}$                      \\ \hline
\end{tabular}
\end{table*}

\begin{algorithm*}
\caption{Algorithm to Directly Calculate the fast Rotation Angles for Non-Symmetric Matrices.}
\label{alg1}
\small
\begin{multicols}{2}
\begin{algorithmic}
\Input {$\mathbf{A}$.}
\Output	{Rotation Matrix $\mathbf{R}_{\tilde{\theta}}$  and $\mathbf{R}_{\tilde{\Theta}}$.}
\Step{1} Calculate the initial values:
$$S_{N1}=Sign(c+b)$$
$$S_{D1}=Sign(d-a)$$
$$N1=|c+b|$$
$$D1=2\times|d-a|$$
$$K1=exp_2(D1)-exp_2(N1)$$
$$S_{N2}=Sign(c-b)$$
$$S_{D2}=Sign(d+a)$$
$$N2=|c+b|$$
$$D2=2\times|d-a|$$
$$K2=exp_2(D2)-exp_2(N2)$$
\Step{2} Calculate $l_{\alpha}$ and $l_{\beta}$ using following case statement:
$$(l1_{temp},B)=\begin{cases}
(K1+1,1) & \text{if } 1.5\times D1> 2^{K1+1}\times N1 \\ 
(K1-1,0) & \text{if } 1.5\times D1< 2^{K1}\times N1 \\ 
(K1,1) & \text{if } D1< 2^{K1}\times N1 \\ 
(K1,0) & \text{default } 
\end{cases}$$
$$l_{\alpha}=max(l1_{temp}+1,2)$$
$$(l2_{temp},b)=\begin{cases}
(K2+1,1) & \text{ if } 1.5\times D2> 2^{K2+1}\times N1 \\ 
(K2-1,0) & \text{ if } 1.5\times D2< 2^{K2}\times N1 \\ 
(K2,1) & \text{ if } D2< 2^{K2}\times N1 \\ 
(K2,0) & \text{ default } 
\end{cases}$$
$$l_{\beta}=max(l2_{temp}+1,2)$$

\Step{3} Calculate the $l_{\theta}$ and $l_{\Theta}$ using the following case statement:
$$(l_{\Theta},l_{\theta})=\begin{cases}
(l_\alpha,l_\alpha) & \text{if } (N2=0)\\ 
(l_\beta,l_\beta) & \text{if } (N1=0)\\ 
(l_\beta-(B\&b),l_\alpha) & \text{if } (l_\beta-l_\alpha=-1)\\ 
(l_\alpha-(B\&b),l_\beta) & \text{if } (l_\beta-l_\alpha=1)\\ 
(l_\beta-1,0) & \text{if } (l_\beta-l_\alpha=0)\\ 
(min(l_\alpha,l_\beta),min(l_\alpha,l_\beta)) & \text{default} 
\end{cases}$$
$$(l_{\Theta},l_{\theta})=\begin{cases}
(l_{\theta},l_{\Theta}) & \text{if } (S_{D2}*S_{N2} \neq S_{D1}*S_{N1})\\ 
(l_{\Theta},l_{\theta}) & \text{default} 
\end{cases}$$
\Step{4} Calculate the $S_{\tilde{{\theta}}}$ and $S_{\tilde{\Theta}}$ using the following case statement:
$$S=\begin{cases}
S_{D1}*S_{N1} & \text{if } (l_\beta-l_\alpha\geqslant 0 \  || \  N2=0)\\ 
S_{D2}*S_{N2} & \text{default} 
\end{cases}$$
$$(S_{\tilde{{\Theta}}},S_{\tilde{{\theta}}})=\begin{cases}
(S,S) & \text{if } (N2=0)\\ 
(S,-S) & \text{if } (N1=0)\\ 
(S*Sign(l_\beta-l_\alpha),S) & \text{default }\\ 
\end{cases}$$
\Step{5} Calculate the $(c_{\tilde{\theta}},s_{\tilde{\theta}})$ and $(c_{\tilde{\Theta}},s_{\tilde{\Theta}})$ pairs using the following case statement:
$$\tilde{t_1}=\begin{cases}
0 & \text{if } l_{\theta}=0\\ 
2^{-{l_{\theta}}} & \text{ default } 
\end{cases}$$
$$\tilde{t_2}=\begin{cases}
0 & \text{if } l_{\Theta}=0\\ 
2^{-l_{\Theta}} & \text{ default } 
\end{cases}$$
$$(c_{\tilde{\theta}},s_{\tilde{\theta}})=\frac{1}{1+\tilde{t_1}^2} \times (1-\tilde{t_1}^2,2 \times S_{\tilde{\theta}} \times \tilde{t_1}) $$
$$(c_{\tilde{\Theta}},s_{\tilde{\Theta}})=\frac{1}{1+\tilde{t_2}^2} \times (1-\tilde{t_2}^2,2 \times S_{\tilde{\Theta}} \times \tilde{t_2}) $$
\Step{6} Calculate the rotation matrices $\mathbf{R}_{\tilde{\theta}}$ and $\mathbf{R}_{\tilde{\Theta}}$ using the following case statement:
$$ \mathbf{R}_{\tilde{\theta}} =\begin{cases}
\begin{bmatrix} c_{\tilde{\theta}} & s_{\tilde{\theta}} \\-s_{\tilde{\theta}} & c_{\tilde{\theta}} \end{bmatrix} & \text{ if } S_{D1}<0\\ \\
\begin{bmatrix} s_{\tilde{\theta}} & c_{\tilde{\theta}} \\c_{\tilde{\theta}} & -s_{\tilde{\theta}} \end{bmatrix} & \text{ else } 
\end{cases}$$
$$ \mathbf{R}_{\tilde{\Theta}} =\begin{cases}
\begin{bmatrix} c_{\tilde{\Theta}} & s_{\tilde{\Theta}} \\-s_{\tilde{\Theta}} & c_{\tilde{\Theta}} \end{bmatrix} & \text{ if } S_{D1}<0\\ \\
\begin{bmatrix} s_{\tilde{\Theta}} & c_{\tilde{\Theta}} \\c_{\tilde{\Theta}} & -s_{\tilde{\Theta}} \end{bmatrix} & \text{ else } 
\end{cases}$$
\end{algorithmic}
\end{multicols}
\end{algorithm*}

\subsection{Reducing the Direct Estimation Complexity}
\label{RDEC}

\begin{figure}
\centering
\includegraphics[width=0.49\textwidth,keepaspectratio]{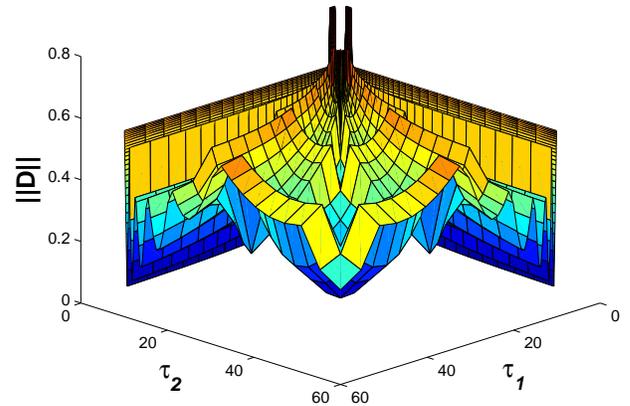}{}
\caption{$\mid\mid D \mid\mid$ vs. $\tau_1$ and $\tau_2$ for ERFHSVD2}
\label{fig:DVsT1T2DirectRelaxedSimp2}
\end{figure}

\begin{figure}
\centering
\includegraphics[width=0.49\textwidth,keepaspectratio]{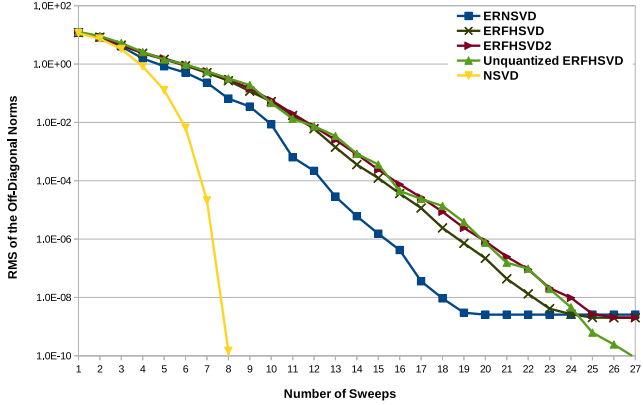}{}
\caption{Comparison of the $RMS_{ODN}$ for NSVD, ERNSVD, ERFHSVD, Unquantized ERFHSVD, and ERFHSVD2}
\label{fig:NormDirectRelaxedSimp}
\end{figure}

Algorithm ~\ref{alg1} shows the complete flow for the calculation of $\tilde{\theta}$ and $\tilde{\Theta}$. The complexity of angle calculation is optimized for achieving reasonable hardware complexity. In \textbf{Step2} the calculation of $(l1_{temp},B)$ or $(l2_{temp},b)$ can be reduced to two comparison if we do not use \emph{B} and \emph{b}, and it will only affect one of the fast rotation angles in the TABLE~\ref{tab:Direct}. This case statement is represented in (\ref{eq:relax1}). We named this method ERFHSVD2. The $\mid\mid D \mid\mid$ of ERFHSVD2 vs. $\tau_1$ and $\tau_2$ is demonstrated in Figure~\ref{fig:DVsT1T2DirectRelaxedSimp2}. The comparison between the  $\mid\mid D \mid\mid$  of ERFHSVD and ERFHSVD2 shows higher values for ERFHSVD which is expected since ERFHSVD2 is the less complex approximation. The $RMS_{ODN}$ of both methods are represented in Figure~\ref{fig:NormDirectRelaxedSimp}. The loss in convergence speed is two extra rotation at maximum performance for 32 bit representation, while it requires an extra comparator.
The performance of unquantized representation of ERFHSVD2 is also demonstrated in Figure~\ref{fig:NormDirectRelaxedSimp} to prove that the floor in the $RMS_{ODN}$ is due to the quantization and not approximations. This figure also demonstrates the difference between ERNSVD and ERFHSVD. The loss in convergence speed is four extra rotations at maximum performance for 32 bit representation, this difference is smaller when larger $RMS_{ODN}$ is acceptable. 
We must note that depending on the implementation, one iteration of ERFHSVD might be equal to applying two iteration of ERNSVD. The ERNSVD needs to calculate and apply the symmetrizing rotation and then calculate and apply the diagonalizing rotation.

\begin{figure}
\centering
\includegraphics[width=0.49\textwidth,keepaspectratio]{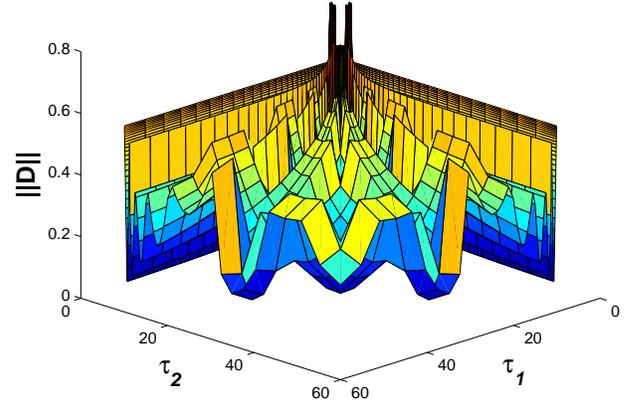}{}
\caption{$\mid\mid D \mid\mid$ vs. $\tau_1$ and $\tau_2$ for ERFHSVD}
\label{fig:DVsT1T2DirectRelaxedSimp}
\end{figure}

\section{Hardware Implementation}
\label{HWI}
\begin{figure}
\centering
\includegraphics[width=0.49\textwidth,keepaspectratio]{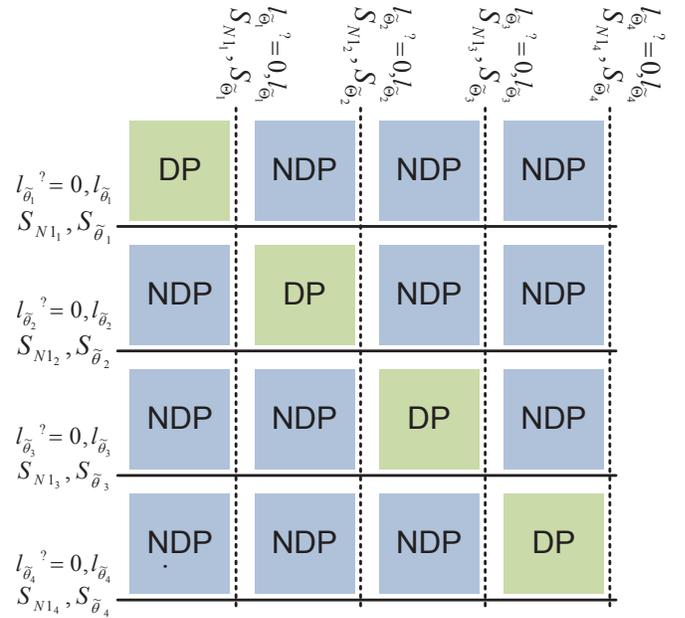}{}
\caption{Architecture of the design.}
\label{fig:Arch}
\end{figure}

These algorithms can be implemented in different methods on the higher level. Figure~\ref{fig:Arch} shows a high level systolic implementation of the decomposition algorithm for an $8 \times 8$ matrix. Figure~\ref{fig:Sched} shows how the scheduling for this implementation can be organized so four independent rotations are applied in each clock cycle. Each pair shows the number of rows and columns that rotation is calculated from and applied to.  While different high level designs choose different method to manage their memory unit, timing, and connections, majority of these architectures follow some common design footsteps. These architectures are constructed from two different types of processor:  diagonal and non-diagonal processors. The diagonal processor calculate, transmit and applies the horizontal and vertical rotations while the non-diagonal processor receives, and applies the rotations. Since all our contributions can be explained with more details in lower design discussions, we focus on the design of basic circuits for diagonal and non-diagonal processors (DP, and NDP). The discussed and presented are only one possible implementation. The two's complement representation is used in this deign to represent negative numbers. The fully combinational implementation of the design is discussed here.

\begin{figure}
\begin{align*}
(p,q)=&(1,2)(3,4)(5,6)(7,8)\\
&(1,4)(2,6)(3,8)(5,7)\\
&(1,6)(4,8)(2,7)(3,5)\\
&(1,8)(6,7)(4,5)(2,3)\\
&(1,7)(8,5)(6,3)(4,2)\\
&(1,5)(7,3)(8,2)(6,4)\\
&(1,3)(5,2)(7,4)(8,6) 
\end{align*}
\caption{Scheduling of matrix processing order.}
\label{fig:Sched}
\end{figure}

\subsection{Calculating the Rotations}
\label{Step1}
\begin{figure}
\centering
\includegraphics[width=0.49\textwidth,keepaspectratio]{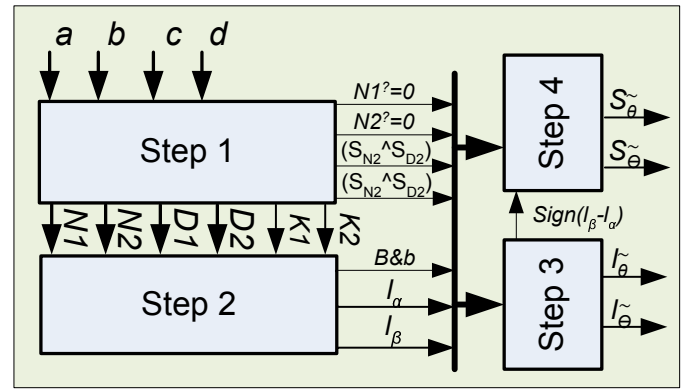}{}
\caption{Diagram for calculating the Given's Rotations based on the ERFHSVD algorithm}
\label{fig:CR}
\end{figure}

Figure~\ref{fig:CR} demonstrates the diagram for calculating the Given's Rotations based on the ERFHSVD algorithm. This involves the first four steps of Algorithm~\ref{alg1}. The inputs are the elements of $2 \times 2$ matrix that is the target of decomposition. The outputs are the sign ($Sign$) and power ($l$) in $tan(\tilde{x})=\tilde{t}=Sign \times2^{-l}$ for both rotation angles $\tilde{\theta}$ and $\tilde{\Theta}$. This circuit is only exists in diagonal processors (DPs). After calculating the value this circuit transmits the values to the circuit for applying rotations in the DP. This circuit also sends the required signals to the circuit for applying rotations in the same DP ($l_{\tilde{\Theta}}$, $l_{\tilde{\theta}}$, $S_{\tilde{\Theta}}$, $S_{\tilde{\theta}}$,$S_{N1}$, $l_{\tilde{\Theta}}^?=0$, and $l_{\tilde{\theta}}^?=0$) and NDPs in the same row ($l_{\tilde{\Theta}}$, $S_{\tilde{\Theta}}$, $S_{N1}$, and $l_{\tilde{\Theta}}^?=0$)\footnote{$l_{\tilde{\Theta}}^?=0$ is the signal name and will be one if $l_{\tilde{\Theta}}$ is equal to zero (this signal is used to make decision on cases in \textbf{Step 2} of the \textbf{Algorithm ~\ref{alg3}}), and $S_{N1}$ is the sign of the numerator in equation (\ref{eq:2}).}, and column ($l_{\tilde{\theta}}$, $S_{\tilde{\theta}}$, $S_{N1}$, and $l_{\tilde{\theta}}^?=0$). The circuit for each of the steps is discussed hereafter.

\subsubsection{ERFHSVD Step 1}
\label{Step1}
\begin{figure}
\centering
\includegraphics[width=0.49\textwidth,keepaspectratio]{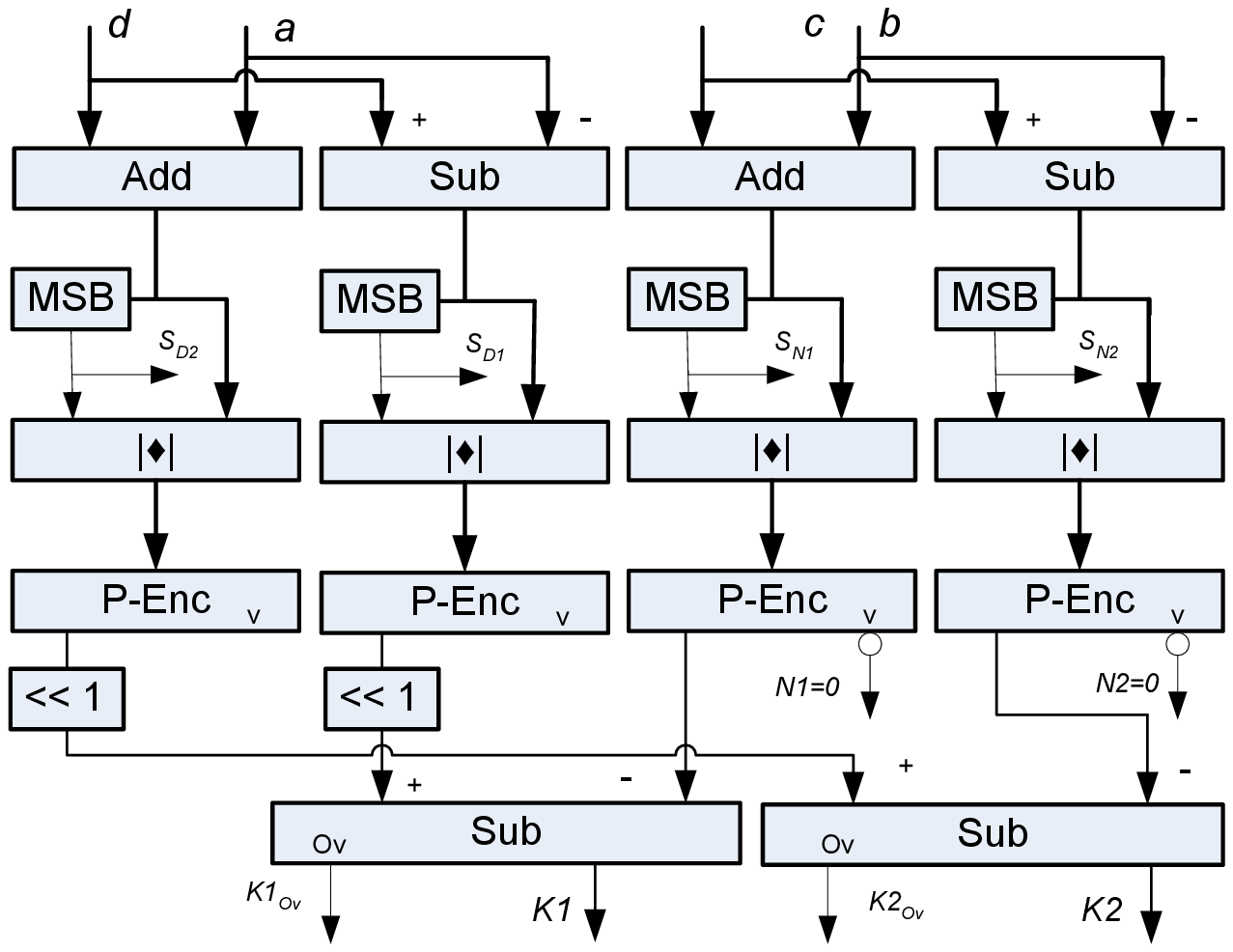}{}
\caption{Diagram of the first step of ERFHSVD algorithm}
\label{fig:Step1}
\end{figure}

Figure~\ref{fig:Step1} demonstrates the diagram of the proposed design for the first step of ERFHSVD algorithm. This circuit generates the initial values for being used at the next steps in ERFHSVD algorithm.  The inputs are the elements of $2 \times 2$ matrix that is the target of decomposition. The MSB block does not have any cost in VLSI implementation and it is demonstrating that the most significant bit of the value is used to determine the sign (The implementation is two's complement). The blocks with $|\blacklozenge|$ on them are the circuits to calculate the absolute value of the input value. We assume this will cost an XOR complimenting circuit and an adder to calculate the two's complement of negative numbers. To avoid overflow or the need to use saturation, the adder size of the $|\blacklozenge|$ circuits should be at least of the same size as the matrix input argument. The P-Enc blocks are priority encoders as explained in (\ref{eq:7}), where $v$ is the valid signal and is zero when the input signal is equal to zero, and it is one for rest of the cases. The blocks with $<<\ 1$ are indicating shifts to the left (or multiplying by two) and their VLSI implementation does not have any cost. The last two subtractors in this diagram have the size of $ceil(log_2(bit))$ where the $bit$ is the number of bits each element of input matrix is represented with. One of the conditions of the second step is to limit the value of $K1$ and $K2$ to the minimum of two. The over flow outputs of the last two subtractors are indicating a negative result since both inputs are positives.\par

\subsubsection{ERFHSVD Step 2}
\label{Step2}
\begin{figure}
\centering
\includegraphics[ angle=90, width=0.49\textwidth, keepaspectratio]{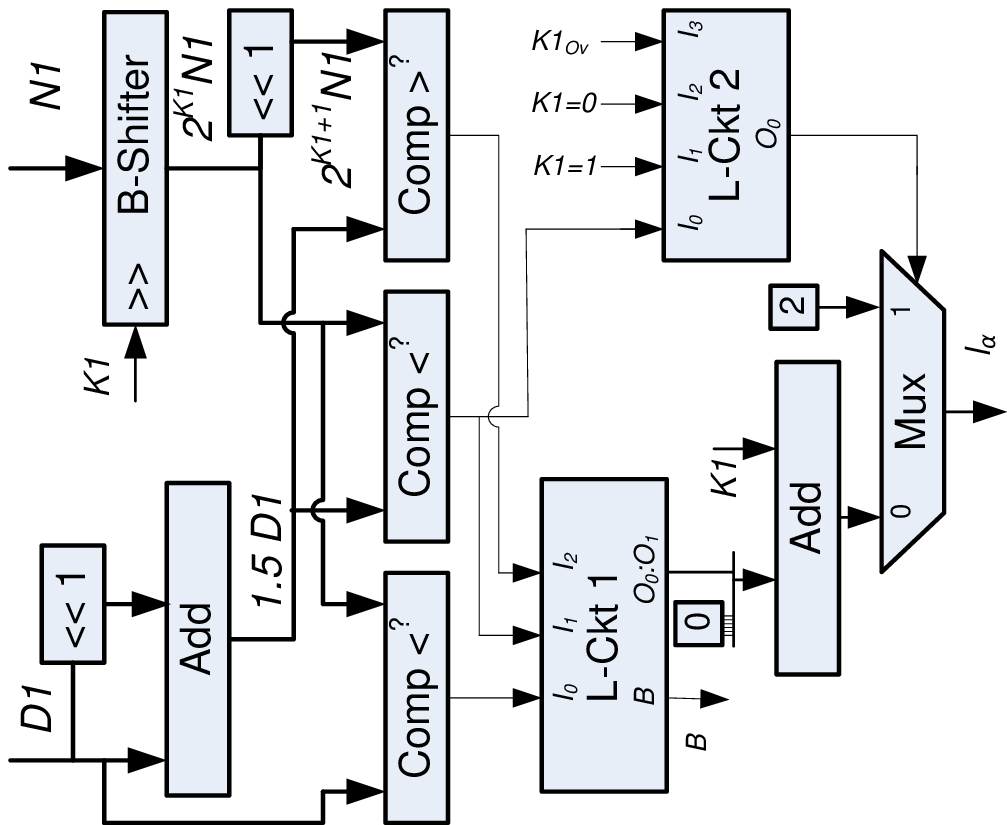}{}
\caption{Diagram of the second step of ERFHSVD algorithm}
\label{fig:Step2}
\end{figure}

Figure~\ref{fig:Step2} demonstrates the diagram of the proposed design for the second step of ERFHSVD algorithm. This block shows how $l_\alpha$ is calculated. A similar block can be used to generate $l_\beta$. This circuit uses an adder and a B-Shifter to generate the signals required for the case statement in \textbf{Step 2}. The blocks with $<<\ 1$ are indicating shifts to the left (or multiplying by two) and their VLSI implementation does not have any cost. The adder, B-Shifter, and comparators are of the size $bit$. The last adder and mux in the block diagram are of the size $ceil(log_2(bit))$. The L-Ckt is indicating a logical circuit that can be implemented with a and, or, and inverter representation or any other means necessarily. We have merged the case statements of \textbf{Step 2} with the mathematics phrase comeing right after them and represented them both in one circuit. In fact the last multiplexer in the diagram is saturating the results to the minimum of two; while the "L-Ckt 2" is generating its select signal inputs. The "L-Ckt 1" outputs, zero, one or two to be added in the final adder based on the case condition; while also generating the "$B$" signal. For "L-Ckt 1", $B= I_0 \overline{I_1}+I_2$, $O_0= \overline{I_1}\ \overline{I_2}$, and $O_1= I_2$.  For "L-Ckt 2", $O_0= I_3 + I_2 + I_1I_0 $. The signal $K1=0$ is generated with an eight-input NOR gate while the signal $K=1$ is generated with a NOR gate and an inverter. We did not show these gates in the figure in favor of keeping the diagrams straightforward.\par

\subsubsection{ERFHSVD Step 3}
\label{Step3}
\begin{figure}
\centering
	\includegraphics[angle=90, width=0.49\textwidth,keepaspectratio]{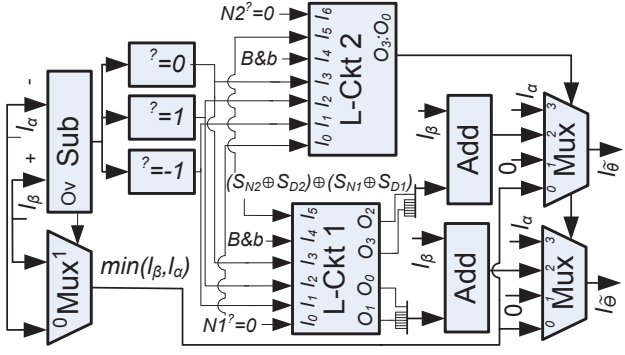}{}
\caption{Diagram of the third step of ERFHSVD algorithm}
\label{fig:Step3}
\end{figure}

Figure~\ref{fig:Step3} demonstrates the diagram of the proposed design for the third step of ERFHSVD algorithm. This circuit takes the $l_\alpha$, $l_\beta$, $N2^?=0$, $B b$, $(S_{N2} \oplus S{D2}) \oplus(S_{N1} \oplus S_{D1})$, $N1^?=0$\footnote{$N1^?=0$ will be one if $N1=0$ (which $N1$ is the numerator in equation(\ref{eq:2})).} signals as inputs and outputs the values of $l_{\tilde{\theta}}$ and $l_{\tilde{\Theta}}$ to the circuit for applying rotations and the sign bit of $l_\beta - l_\alpha$ to \textbf{Step~4}. This design  merges the case statement in \textbf{Step~3} and the mathematical phrase after it, and implements both together. Shift left logical (SLL blocks marked with $<<$), adder, subtractor, and multiplexers blocks are of the size $ceil(log_2(bit))$. The block that is supposed to determine if the result of the subtract is zero (this block is marked with $^?=0$) requires a NOR gate of size $ceil(log_2(bit))+1$ assuming that overflow (Ov) can happen and for detecting "one" an inverter and a NOR gate is required, since both $l_\alpha$ and $l_\beta$ are positive an AND gate of size $ceil(log_2(bit))$ with an inverter can be used to synthesize the block marked with $^?=-1$. This provision is taken to prevent the usage of comparators. The L-Ckt~1 applies any change necessary on the value of $l_\beta$ by adding zero, minus one, or minus two to it. The L-Ckt~2 generates the input signals to the mux based on the inputs to assure that correct values are assigned to $l_{\tilde{\theta}}$, and $l_{\tilde{\Theta}}$. Since $l_\beta \geq 2$, the result of the additions can not be negative, the Ov outputs of the adders can be ignored. The final multiplexer of the circuit is needed to guarantee that diagonalized output matrix is normalized. In L-Ckt~1, $O_0= \overline{I_0} \ \overline{I_1} \ \overline{I_2} \ \overline{I_5}+\overline{I_0} \ I_2 \ \overline{I_4} \ \overline{I_5}+\overline{I_0} \ I_1 \ I_4 \ \overline{I_5}$, $O_1= \overline{I_1} \ I_4 \ \overline{I_5}+\overline{I_0}\ \overline{I_1} \ \overline{I_5}$, $O_2= \overline{I_0} \ \overline{I_1} \ \overline{I_2} \ I_5+\overline{I_0} \ I_2 \ \overline{I_4} \ I_5+\overline{I_0} \ I_1 \ I_4 \ I_5$, and $O_3= \overline{I_1} \ I_4 \ I_5+\overline{I_0}\ \overline{I_1} \ I_5$. In L-Ckt~2, the signals are defined as follows: $O_0= I_5 \ I_3 + I_5 \ \overline{I_0} \ I_1 + I_6$, $O_1= I_0+I_1+I_2+I_3 \overline{I_5}+I_6$, $O_2= I_3 \overline{I_5}+  \overline{I_0}\ I_1 \ \overline{I_5} + I_6$, and $O_3= I_0+I_1+I_2+I_5 \ I_3+ \ I_6$.\par

Using the circuit presented in Figure~\ref{fig:Step3} with minor changes in L-Ckt~1 and L-Ckt~2 implementation, every different value in the first three rows of TABLE~\ref{tab:Direct} can be assigned as the rotation angle. If implementation of more rows from the table is required, the circuit to detect plus and minus two also should be added and fed to the logical circuits. The simple design presented for this step of the algorithm would assist the designers to change table values based on their application need. For more complex design, in case the number of inputs to the logical circuits is high, an alternative design can be explored that two-input multiplexers (two AND gates of size $ceil(log_2(bit))$) are used to impediment the second case statement in \textbf{Step~3}. The multiplexer swaps the value of $l_{\tilde{\theta}}$, and $l_{\tilde{\Theta}}$ if $(S_{D1}\oplus S_{N1}) \oplus (S_{D2} \oplus S_{N2})=1$. This will reduce the complexity of L-Ckt~1 and might eliminate the need for one of the adders.\par 

\subsubsection{ERFHSVD Step 4}
\label{Step4}
\begin{figure}
\centering
	\includegraphics[angle=90, width=0.3\textwidth,keepaspectratio]{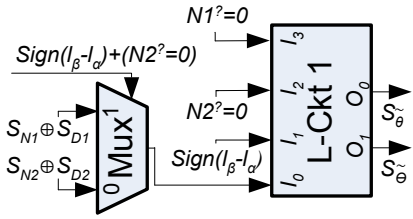}{}
\caption{Diagram of the fourth step of ERFHSVD algorithm}
\label{fig:Step4}
\end{figure}

Figure~\ref{fig:Step4} demonstrates the diagram of the proposed design for the fourth step of ERFHSVD algorithm. This circuit takes the sign bit of $l_\beta - l_\alpha$ ($Sign(l_\beta - l_\alpha)$), $N2^?=0$, $N1^?=0$, $S_{N1}$, $S_{D1}$, $S_{N2}$, and $S_{D2}$ as input and generates the sign of the rotations ($S_{tilde{\Theta}}$ and $S_{tilde{\theta}}$). Correct implementation of this step is vital for convergence of the iterations as well as achieving the normalized diagonal elements. For "L-Ckt~1", $O_0= \overline{I_0} \  \overline{I_2} I_3 + I_0 \overline{I_2} \ \overline{I_3}+ I_0 I_2$, and $O_1= \overline{I_0} I_1 \overline{I_2}\ \overline{I3}+I_0I_3+ I_0I_2+I_0 \overline{I_1} $.

\subsection{Applying the Rotations}
\label{AR}
\begin{figure}
\centering
\includegraphics[keepaspectratio,width=0.49\textwidth]{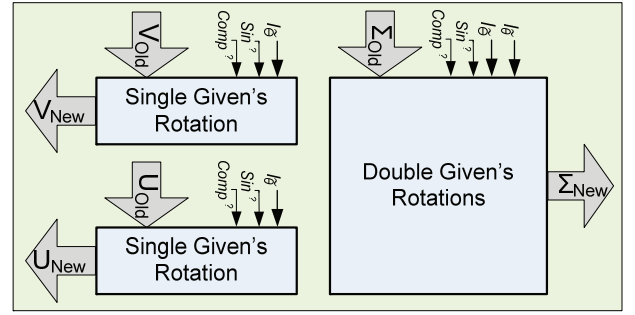}{}
\caption{Diagram of the circuit for applying rotations}
\label{fig:AR}
\end{figure}

Figure~\ref{fig:AR} demonstrates the diagram of the proposed design for applying rotations. This circuit includes the major part of the NDP (except memory units and the mechanism to receive the input arguments and transfer the results). Assuming a $2 \times 2$ matrix decomposition is presented in (\ref{eq:Decomp}). For an iterative algorithm this unit should be able to apply both rotation matrices on the old value of $\pmb{\Sigma}$ ($\pmb{\Sigma}_{Old}$) based on the values of $l_{\tilde{\theta}}$ and $l_{\tilde{\Theta}}$, it also should apply one rotation matrix on the old value of $\mathbf{U}$ ($\mathbf{U}_{Old}$) based on the values of $l_{\tilde{\theta}}$ and one rotation matrix on the old value of $\mathbf{V}$ ($\mathbf{V}_{Old}$) based on the values of $l_{\tilde{\Theta}}$ to generate the new values. The $Sin^?$\footnote{$Sin^?$ is one when a block going to apply sine rotation and cosine rotation otherwise.} signal is in fact one if $S_N \geq 0$ as we will explain in the next sub section. The initial value of $\pmb{\Sigma}$ is $\mathbf{A}$ and the initial value of $\mathbf{U}$ and $\mathbf{V}$ is identity matrix. \par

\begin{equation}
\mathbf{A}= \mathbf{U} \times \pmb{\Sigma} \times \mathbf{V}^T 
\label{eq:Decomp}
\end{equation}

\subsubsection{Applying Double and Single Given's Matrix Rotations}
\label{M&S}
\begin{figure}
\centering
\includegraphics[width=0.49\textwidth,keepaspectratio]{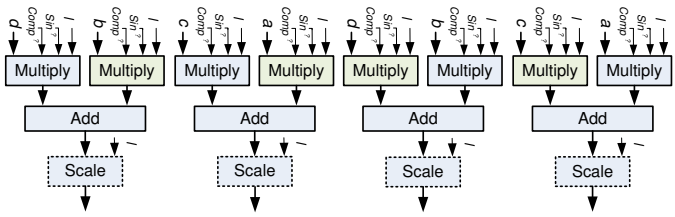}{}
\caption{Diagram for applying single Given's rotation}
\label{fig:SGR}
\end{figure}

Figure~\ref{fig:SGR} demonstrates the diagram of the proposed design for calculating the single Given's Rotation. It is important to note that applying a Given's rotation is in fact a multiplication of two $2 \times 2$ matrices. This is equal to eight multiplication and four addition if no other consideration is made about the implementation of the system. We can use this circuit to calculate the value of $V$ and $U$. The scale circuit is demonstrated with doted outline to remind readers that its presence or implementation complexity (and accuracy as the result) can be adapted based on the application need. The $Sin^?$ signal that is an input to each Multiply block shows if that block has to multiply the input argument with $Sin(\tilde{x})$ or $Cos(\tilde{x})$. A control unit can separately assign different values to each Multiply block, or it can choose to multiply the input value with $Sin(\tilde{x})$ or $Cos(\tilde{x})$. In the design the blocks are with a different color to make sure they multiply the input argument with $Cos(\tilde{x})$ while the normal blocks are multiplying the input arguments with $Sin(\tilde{x})$. \par

\begin{figure*}
\centering
\includegraphics[width=0.95\textwidth,keepaspectratio]{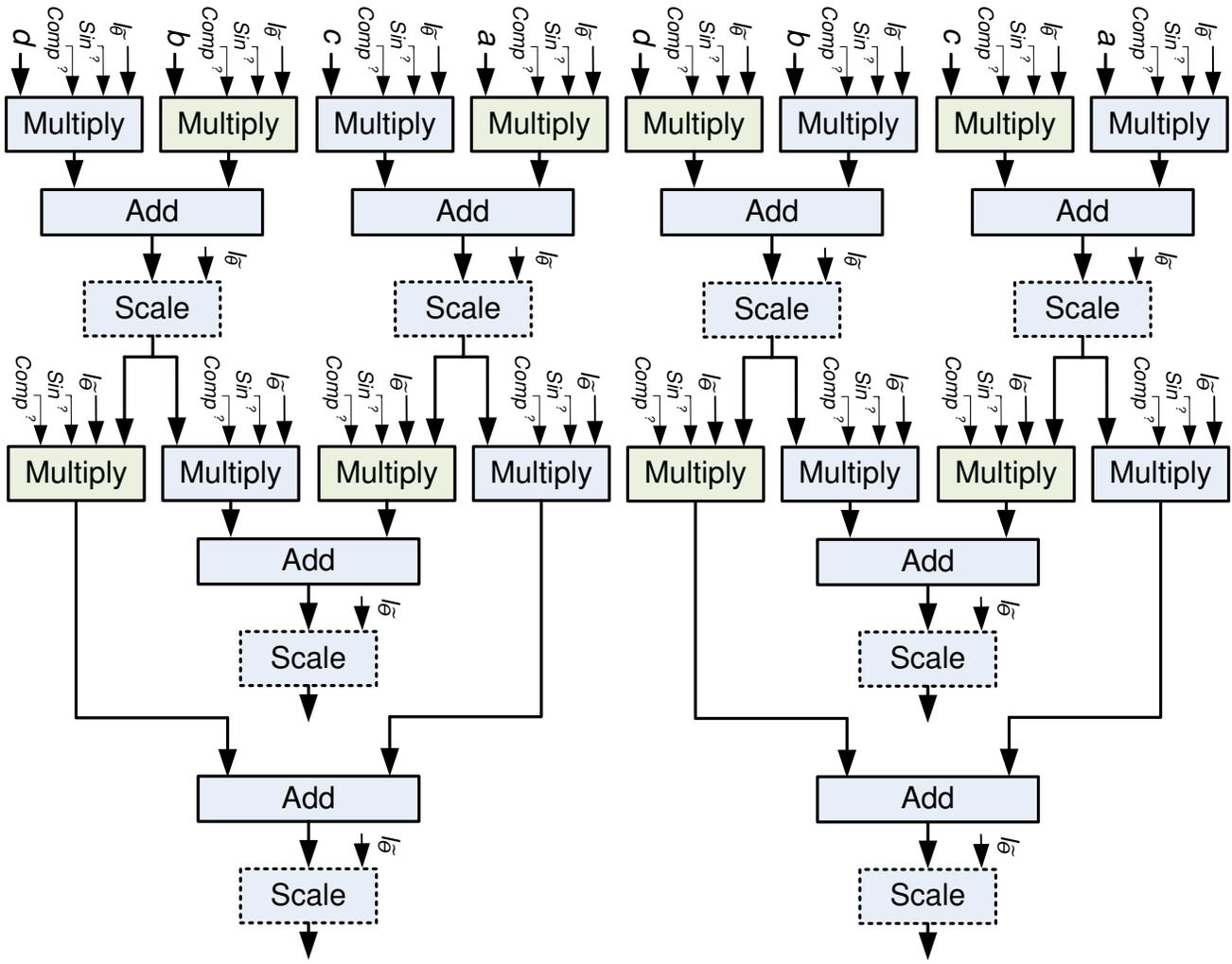}{}
\caption{Diagram for applying double Given's rotations}
\label{fig:DGR}
\end{figure*}

Figure~\ref{fig:DGR} demonstrates the diagram of the proposed design for calculating the single Given's Rotations. The same notation as in Fig~\ref{fig:SGR} is used in this diagram. The circuit in this diagram has double complexity compared to the design of Figure~\ref{fig:SGR}. This is expected since this circuit has to multiply the input $2 \times 2$ matrix with two rotation matrices of $R_{\tilde{\theta}}$ and $R_{\tilde{\Theta}}$. \par

\subsubsection{Multiplying and Scaling}
\label{M&S}
The rotation angles are the approximations to $t=|tan(x)|$ with $\tilde{t}=|tan(\tilde{x})|=2^{-l}$ and as the result of applying double rotations the elements of any rotation matrix should be zero, one, or the values derived from equation (\ref{eq:cos}) or equation (\ref{eq:sin}). The $Sign$ in (\ref{eq:sin}) will be determined in \textbf{Step 4} ($S_{\tilde{\theta}}, \ or \ S_{\tilde{\Theta}}$). The part $\frac{1}{1+\tilde{t}^2}$ is common between all the coefficients and we will discuss it in detail after presenting the circuit which applies the uncommon part of the multiplications (this is refereed to as scaling) if the proposed circuit is able to apply any of these coefficients using only control signals then we can merge \textbf{Step 5} and \textbf{Step 6} in algorithm \ref{alg1} (this is refereed to as multiplying).\par

\begin{equation}
\frac{1}{1+\tilde{t}^2} \times 1-\tilde{t}^2
\label{eq:cos}
\end{equation}

\begin{equation}
\frac{1}{1+\tilde{t}^2} 2 \times Sign \times \tilde{t}
\label{eq:sin}
\end{equation}

\begin{figure}
\centering
\includegraphics[angle=90,width=0.38\textwidth,keepaspectratio]{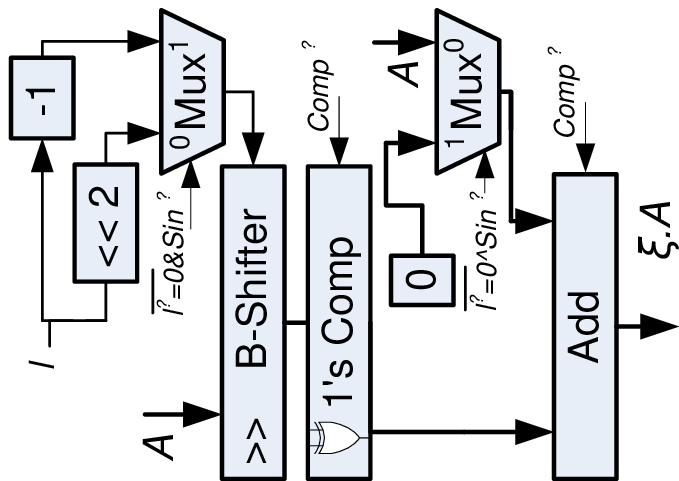}{}
\caption{Diagram of the multiplier for ERFHSVD algorithm}
\label{fig:Multiply}
\end{figure}

Figure~\ref{fig:Multiply} demonstrates the diagram of the proposed design for multiplications. This is a fully combinational circuit and a possible design. The input $A$ is any of the elements of the matrix in (\ref{eq:1}). If $\xi$ is zero ($l=0, \ Sin^?=1, \  Comp^?=1$), then the adder in the figure is adding $A$ with $-A$. If $\xi$ is one ($l=0, \ Sin^?=0, \  Comp^?=0$), then the adder is adding $A$ with zero. To apply equation (\ref{eq:cos}), the input value has to be shifted to the right twice the value of $l$ and then subtracted from the original value ($l \neq 0, \ Sin^?=0, \  Comp^?=1$). The circuit should be able to apply equation (\ref{eq:sin}) assuming $sign$ can be negative or positive. For this equation the adder only adds one to the input value to convert one's compliment values to two's compliment values ($l \neq 0, \ Sin^?=1, \  Comp^?=1$). For this circuit a control unit can define the values for selecting input of multiplexers or the value of $Comp^?$\footnote{$Comp^?$ is one if the multiplicand is negative.}; however if we use the signals as defined in Figure~\ref{fig:Multiply}, we can assign $Sin^?=S_{N1}$ and just swap the input index for multiplexers if a unit needs to apply $Cos(\tilde{x})$ when $S_{N1} < 0$. If the multiplication circuit belongs to a DP then the $S_{N1}$ signal of the same processor is used to determine the $Sin^?$ signals, while if the circuit belongs to an NDP then $S_{N1}$ that is transmitted from the DP in the same row is used for determining the $Sin^?$ signals in the first row of the double Given's rotations circuit and updating the values of $\mathbf{U}$. The $S_{N1}$ that is transmitted from the DP in the same column is used for determining the $Sin^?$ signals in the second row of the double Given's rotations circuit and updating the values of $\mathbf{V}$.\par

\begin{figure}
\centering
\includegraphics[ width=0.35\textwidth,keepaspectratio]{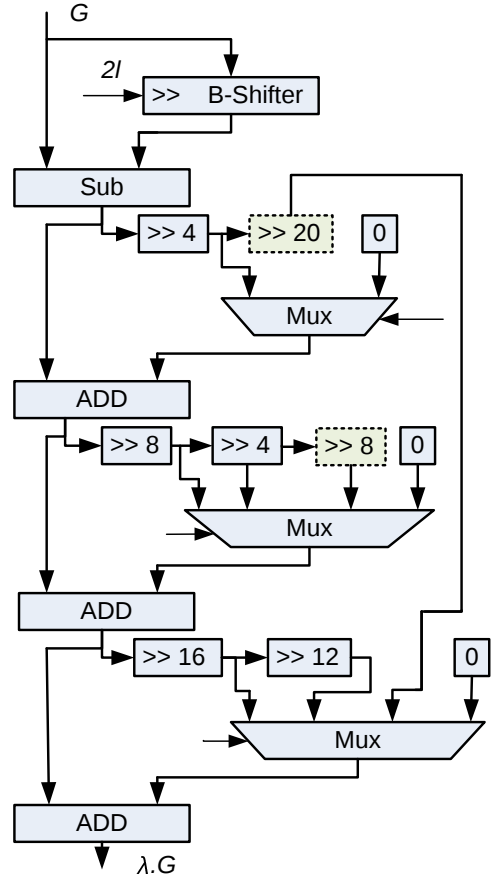}{}
\caption{Diagram of the scaling circuit for ERFHSVD algorithm}
\label{fig:Scaling}
\end{figure}

Applying scaling is one of the more resource demanding parts in this design while it might not be necessary for every application to apply scaling coefficients when decomposing a matrix. This generally results in applying orthogonal rotation angles and an increase in the value of diagional elements of decomposition which might be acceptable for some applications. Comparing the fast rotations method with CORDIC method the complexity of the scaling circuit is increased when applying fast rotations. The fact that CORDIC method applies all the rotation angles and only the sign of the rotations are different results in a constant scaling value; while, the fast rotations can be different any time, this results in different scaling for different rotation angels. \cite{gotze1993efficient} and \cite{gotze1991parallel} suggest using the Tailor series representation of $\frac{1}{1+\tilde{t}^2}$ as in (\ref{eq:scale}). The estimation of complexity depends on the implementation, but a simplified look at the problem is presented in \cite{gotze1993efficient} and \cite{gotze1991parallel}. The complexity of other CORDIC based implementations is also presented there. The implementation of the scaling factor for 32 bits implementation requires four Shifts and four additions. A closer look at the scaling coefficient would help reducing the complexity of scaling circuit to four addition, and one shift. Table~\ref{tab:scaling} shows the scaling coefficient that needs to be applied for each rotation angle of $tan(\tilde{x})=2^{-l}$. Increase in the value of $l$ causes the scaling coefficient to merge to one; as the result, for $l \geq 16$ the 32 bit representation of the coefficient is rounded to one. The notation $Z= Acc.\{\Delta\} \gg \delta$ is equal to $Z = \Delta+\Delta>>\delta$ which $\Delta \gg \delta$ means arithmetic shift of the $\Delta$ to right $\delta$ times.  
Figure~\ref{fig:Scaling} demonstrates the diagram of the proposed design for scaling in ERFHSVD algorithm.\par

\begin{equation}
\frac{1}{1+\tilde{t}^2}=(1-2^{-2l})(1+2^{-4l})(1+2^{-8l})(1+2^{-16l})...
\label{eq:scale}
\end{equation}

\begin{table*}[]
\centering
\caption{Scaling values for different rotation angles}
\label{tab:scaling}
\footnotesize
\begin{tabular}{|c|c|c|c|c|c|}
\hline
\multicolumn{5}{|c|}{Scale ($\lambda$) }                                                                                           & \multirow{2}{*}{\begin{tabular}[c]{@{}c@{}}Circuit representation of\\ $\lambda . G $ (any 32 bit input)\end{tabular}} \\ \cline{1-5}
l  & -2l & Fraction                           & Decimal                         & Binary (32 bit)                     &                                                                                                          \\ \hline
\hline
1  & 2   & $\frac{1}{1+\frac{1}{4}}$          & $\frac{4}{5}$                   & 0.11001100110011001100110011001100  & $Acc.\{Acc.\{Acc.\{G-2^{-2}G\}>>4\}>>8\}>>16$                                                           \\ \hline
2  & 4   & $\frac{1}{1+\frac{1}{16}}$         & $\frac{16}{17}$                 & 0.11110000111100001111000011110000  & $Acc.\{Acc.\{G-2^{-4}G\}>>8\}>>16$                                                                      \\ \hline
3  & 6   & $\frac{1}{1+\frac{1}{64}}$         & $\frac{64}{65}$                 & 0.11111100000011111100000011111100  & $Acc.\{G-2^{-6}G\}>>12+(G-2^{-6}G)>>24$                                                                  \\ \hline
4  & 8   & $\frac{1}{1+\frac{1}{256}}$        & $\frac{256}{257}$               & 0.11111111000000001111111100000000  & $Acc.\{G-2^{-8}G\}>>16$                                                                                 \\ \hline
5  & 10  & $\frac{1}{1+\frac{1}{1024}}$       & $\frac{1024}{1025}$             & 0.11111111110000000000111111111100  & $Acc.\{G-2^{-10}G\}>>20$                                                                                \\ \hline
6  & 12  & $\frac{1}{1+\frac{1}{4096}}$       & $\frac{4096}{4097}$             & 0.11111111111100000000000011111111  & $Acc.\{G-2^{-12}G\}>>24$                                                                                \\ \hline
7  & 14  & $\frac{1}{1+\frac{1}{16384}}$      & $\frac{16384}{16385}$           & 0.11111111111111000000000000001111  & $Acc.\{G-2^{-14}G\}>>28$                                                                                \\ \hline
8  & 16  & $\frac{1}{1+\frac{1}{65536}}$      & $\frac{65536}{65537}$           & 0.11111111111111110000000000000000  & $G-2^{-16}G$                                                                                             \\ \hline
9  & 18  & $\frac{1}{1+\frac{1}{262144}}$     & $\frac{262144}{262145}$         & 0.11111111111111111100000000000000  & $G-2^{-18}G$                                                                                             \\ \hline
10 & 20  & $\frac{1}{1+\frac{1}{1048576}}$    & $\frac{1048576}{1048577}$       & 0.11111111111111111111000000000000  & $G-2^{-20}G$                                                                                             \\ \hline
11 & 22  & $\frac{1}{1+\frac{1}{4194304}}$    & $\frac{4194304}{4194305}$       & 0.11111111111111111111110000000000  & $G-2^{-22}G$                                                                                             \\ \hline
12 & 24  & $\frac{1}{1+\frac{1}{16777216}}$   & $\frac{16777216}{16777217}$     & 0.11111111111111111111111100000000  & $G-2^{-24}G$                                                                                             \\ \hline
13 & 26  & $\frac{1}{1+\frac{1}{67108864}}$   & $\frac{67108864}{67108865}$     & 0.11111111111111111111111111000000  & $G-2^{-26}G$                                                                                             \\ \hline
14 & 28  & $\frac{1}{1+\frac{1}{268435456}}$  & $\frac{268435456}{268435457}$   & 0.11111111111111111111111111110000  & $G-2^{-28}G$                                                                                             \\ \hline
15 & 30  & $\frac{1}{1+\frac{1}{1073741824}}$ & $\frac{1073741824}{1073741825}$ & 0.11111111111111111111111111111100  & $G-2^{-30}G$                                                                                             \\ \hline
16 & 32  & $\frac{1}{1+\frac{1}{4294967296}}$ & $\frac{4294967296}{4294967297}$ & 0.111111111111111111111111111111111 & $G-2^{-32}G \approx G$                                                                                           \\ \hline
\end{tabular}
\end{table*}

In Figure~\ref{fig:Scaling} a 32 bit value is multiplied by scale factor $\lambda$ based on the representation in the last column of table~\ref{tab:scaling}. This circuit tries to gain the benefit from repetitive nature of coefficient $\lambda$. If $l \geq 16$ then the scaling factor would simply be one. The circuit presented in this figure is not extendable directly for higher number of bits accuracy and it would be a simpler circuit for lower number of bits. The select signals for the multiplexers has to be assigned based on the value of $l$. For the multiplexer we assumed that the three and four-input multiplexers are made of two-input multiplexers and rearrange them in a way that minimum number of two-input multiplexers are required. In addition, we assume that the Barrel shifters are made out of two-input multiplexers to be able to compare the complexity of the new implementation with the implementation that requires four shifts and four adds. Following the assumptions made, the shifts are requiring the use of B-Shifter since the value of $l$ is different each time. Each B-Shifter will require 160 two-input multiplexers, and our design will require 224 two-input multiplexers for the multiplexers represented in the circuit. This will result in saving 156 two-input multiplexers. This results is approximately saving the cost of one B-Shifter. This circuit also has effects on critical path delay. While a B-Shifter that uses 160 two-input multiplexer has five level of two-input multiplexers the multiplexers in our design at most have two levels. The major benefit of the circuit presented in Figure~\ref{fig:Scaling} is when some error occurs in applying the scaling factor is acceptable. If we implement the circuit only with the first subtraction, the maximum error in applying the scale coefficient will be $6.25\%$. If we represent the circuit with the second adder the maximum error in scaling coefficient will be $0.39\%$ (it multiplies the input value with 0.797 instead of 0.8). The circuit with three adders would have maximum of $0.024\%$ error in applying the coefficients (it multiplies the input value with 0.7998 instead of 0.8). Due to the complexity of scaling circuit in comparison with the rest of the architecture and its variable importance due to the application, we recommend a detailed study for each application. Based on the required accuracy of the result the scaling circuit may be completely ignored or represented with lees complex and less accurate circuit. For example the two shifts that are represented in different color (and doted exterior line), their associated paths, and the last level of adder/multiplexer are only need to be implemented if 32 bit accuracy for the scaling coefficient is needed.\par

\section{Complexity}
\label{Complexity}
The calculation of computational complexity (delay, and resource requirement) depends on the multiple factors other than only the block diagram of the subsystem; however we will try to provide an approximation using some basic assumption to provide a better understanding of this work achievements. We divide the discussions to two sub sections of calculating and applying the rotations. We assume the critical path of any of the blocks in the block diagram representation of any circuit is known as it is represented with symbol $\Delta$. We assume that every subtractor (Sub) is an adder with carry-in equal to one and an array of inverters; as the result, its delay is higher than an adder of same number of bits. We assume that the size of the inputs are $\Lambda$ bits with $\lambda=ceil(log_2\Lambda)$.\par
For resource estimation, we use $A$ as area or resources required for implementing of a particular circuit. The $A_{\Lambda bit Add}$ is representing the area or resources required to implement an adder of size $\Lambda$ bits.

\subsection{Calculating the rotations}
\label{CxCR}
In order to be able to evaluate the delay of the system we assume that an adder of size $\Lambda$ bits has higher delay than a B-Shifter of the same size.

\subsubsection{Step 1}
\label{CxStep1}
The critical path starts from the subtractor that generates the absolute value of $D1$ or $D2$ in the first step and the adder that generates $1.5 \times D1$ in the second stage, instead of the adder that generates the $N1$ or $N2$ in the first Stage and the B-Shifter in the second stage. In (\ref{eq:DStep1}) the critical path delay for the first step of ERFHSVD algorithm is presented.

\begin{equation}
\begin{aligned}
\Delta_{Step~1}=\Delta_{\Lambda bit Sub}+\Delta_{\Lambda bit |\blacklozenge|}+\Delta_{\Lambda bit P-Enc} \\ +\Delta_{\lambda bit Sub}
\end{aligned}
\label{eq:DStep1}
\end{equation}

The phrase $\Delta_{\Lambda bit Sub}$ means that the critical path of a $\Lambda$ bit subtractor; similar notations are used for other elements. The $|\blacklozenge|$ is the circuit calculating the absolute values, and $\Delta_{\Lambda bit P-Enc}$ refers to the critical path of a $\Lambda$ bit P-Enc.
As explained in the beginning of this section the delay path of the subtractor and $|\blacklozenge|$ can be translated to the delay of adders. In (\ref{eq:DStep1-2}) we tried to apply this.\par 

\begin{equation}
\begin{aligned}
\Delta_{Step~1}=\Delta_{\Lambda bit Add}+\Delta_{Neg}+\Delta_{\Lambda bit Add} \\+ \Delta_{Compelement}+ \Delta_{\Lambda bit P-Enc}+\Delta_{\lambda bit Add}+\Delta_{Neg} \\
 = 2 \times \Delta_{\Lambda bit Add}+ 2 \times \Delta_{Inv}+ \Delta_{2 Inp XOR}\\ + \Delta_{\Lambda bit P-Enc}+\Delta_{\lambda bit Add}
\end{aligned}
\label{eq:DStep1-2}
\end{equation}

We assume that $\Delta_{Neg}$ refers to the delay of an inverter and $\Delta_{Compelement}$ refers to the delay of a two-input XOR gate. We replaced the  $\Delta_{Neg}$ with $\Delta_{Inv}$ (inverter) and $\Delta_{Compelement}$ with $\Delta_{2 Inp XOR}$ (two-input XOR) and we will keep that notation hereafter. The Area required to implement this circuit is represented in (\ref{eq:AStep1}).

\begin{equation}
\begin{aligned}
A_{Step~1}=8 \times A_{\Lambda bit Add}+ 2 \times (\Lambda+\lambda) \times A_{Inv} 
\\+4 \times \Lambda \times A_{2 Inp XOR} + 4 \times \Delta_{\Lambda bit P-Enc}
\\+2 \times A_{\lambda bit Add}
\end{aligned}
\label{eq:AStep1}
\end{equation}

\subsubsection{Step 2}
\label{CxStep2}
In (\ref{eq:DStep2}) the Critical path delay for the first step of ERFHSVD algorithm is presented. In this equation $\Delta_{2 Inp Mux}$ is the critical path delay of a two-input multiplexer.\par 

\begin{equation}
\begin{aligned}
\Delta_{Step~2}=\Delta_{\Lambda bit ADD}+\Delta_{Complement}+\Delta_{L-Ckt~1} \\ +\Delta_{\lambda bit Add}+ \Delta_{2 Inp Mux} \\= \Delta_{\Lambda bit ADD}+\Delta_{2 Inp XOR}+\Delta_{L-Ckt~1} \\ +\Delta_{\lambda bit Add}+ \Delta_{2 Inp Mux} 
\end{aligned}
\label{eq:DStep2}
\end{equation}

Where,

\begin{equation}
\begin{aligned}
\Delta_{L-Ckt~1}=\Delta_{Inv}+\Delta_{2 Inp And}+\Delta_{2 Inp Or}
\end{aligned}
\label{eq:DStep2}
\end{equation}

The Area required to implement this circuit is represented in (\ref{eq:AStep2}).

\begin{equation}
\begin{aligned}
A_{Step~2}=2 \times  (A_{\Lambda bit Add}+ A_{\Lambda bit Shifter} \\  + 3 \times A_{\Lambda bit Comparetor} + A_{L-Ckt~1}  + 4 A_{L-Ckt~2}+ \\ A_{\lambda bit Add} + \lambda \times A_{2 Inp Mux})
\end{aligned}
\label{eq:AStep2}
\end{equation}

Where,

\begin{equation}
\begin{aligned}
A_{L-Ckt~1}=3 \times  (A_{Inv}+ A_{2 Inp Or}+ 2 \times A_{2 Inp And}
\end{aligned}
\label{eq:AStep2}
\end{equation}

and

\begin{equation}
\begin{aligned}
A_{L-Ckt~2}= A_{3 Inp Or}+ A_{2 Inp And}
\end{aligned}
\label{eq:AStep2}
\end{equation}

\subsubsection{Step 3}
\label{CxStep3}
In (\ref{eq:DStep3}) the Critical path delay for the first step of ERFHSVD algorithm is presented. In this equation $\Delta_{2 Inp Mux}$ is the critical path delay of a two-input multiplexer.\par 

\begin{equation}
\begin{aligned}
\Delta_{Step~3}=\Delta_{\lambda bit Sub}+\Delta_{\lambda + 1  Inp Nor}+\Delta_{L-Ckt~1} \\ + \Delta_{\lambda bit Add}+ \Delta_{4 Inp Mux} \\= 2 \times \Delta_{\lambda bit Add}+\Delta_{\lambda + 1  Inp Nor}+\Delta_{L-Ckt~1} \\ + \Delta_{Inv}+ \Delta_{4 Inp Mux}
\end{aligned}
\label{eq:DStep3}
\end{equation}

If we assume the four-input multiplexer is made of two levels of two-input multiplexers, we can further simplify the equation.\par 

\begin{equation}
\begin{aligned}
\Delta_{Step~3}= 2 \times \Delta_{\lambda bit Add}+\Delta_{\lambda + 1  Inp Nor}+\Delta_{L-Ckt~1} \\ + \Delta_{Inv}+ 2 \times \Delta_{2 Inp Mux}
\end{aligned}
\label{eq:DStep3-2}
\end{equation}

where,

\begin{equation}
\begin{aligned}
\Delta_{Ckt~1}= \Delta_{Inv} + \Delta_{3 Inp Or}+ \Delta_{4 Inp And}
\end{aligned}
\label{eq:DStep3}
\end{equation}

The Area required to implement this circuit is represented in (\ref{eq:AStep3}). We assume that nine-input NOR and AND gates do require the same area.

\begin{equation}
\begin{aligned}
A_{Step~3}= 3 \times A_{\lambda bit Add}+3\ times A_{(\lambda + 1) Inp Nor}\\ + A_{L-Ckt~1}  
+ A_{L-Ckt~2} + 10 \times A_{Inv}\\ + \lambda \times A_{4 Inp Mux}+ \lambda \times A_{2 Inp Mux}
\end{aligned}
\label{eq:AStep3}
\end{equation}

If we assume that a four-input multiplexer, is made of three, two-input multiplexers the area requirement is presented in (\ref{eq:AStep3-2}).

\begin{equation}
\begin{aligned}
A_{Step~3}= 3 \times A_{\lambda bit Add}+3 \times A_{\lambda + 1  Inp Nor}\\ + A_{L-Ckt~1}  
+ A_{L-Ckt~2} + 10 \times A_{Inv}\\ + \Lambda \times A_{2 Inp Mux}
\end{aligned}
\label{eq:AStep3-2}
\end{equation}

Where,

\begin{equation}
\begin{aligned}
A_{Ckt~1}= 6 \times A_{Inv} + 2 \times A_{3 Inp Or}+ 2 \times A_{2 Inp Or}\\
+ 4 \times A_{3 Inp And} +  6 \times A_{4 Inp And}
\end{aligned}
\label{eq:AStep3-2}
\end{equation}

and,

\begin{equation}
\begin{aligned}
A_{Ckt~2}= 2 \times A_{Inv} + 2 \times A_{3 Inp Or} + 2 \times A_{5 Inp Or}\\
+ 3 \times A_{2 Inp And} +  2 \times A_{3 Inp And}
\end{aligned}
\label{eq:AStep3-2}
\end{equation}

\subsubsection{Step 4}
\label{CxStep4}
In (\ref{eq:AStep4}) the area requirement of the proposed circuit for the fourth step of ERFHSVD is estimated. The Critical path delay of this circuit is not presented since the delay of that circuit does not have any affect on the total delay of the design. 
 In fact the delay of this step is less than the delay of the third step of the algorithm which runs in parts parallel to this step of the ERFHSVD algorithm. \par 

\begin{equation}
\begin{aligned}
A_{Step~4}=  A_{2 Inp Mux}+ A_{L-Ckt~1}+  A_{2 Inp XOR}
\end{aligned}
\label{eq:AStep4}
\end{equation}

where,

\begin{equation}
\begin{aligned}
A_{Ckt~1}= times A_{4 Inp Or} + imes A_{3 Inp Or}\\ 
+ 3 \times A_{3 Inp And} +  4 \times A_{2 Inp And}
\end{aligned}
\label{eq:AStep3-2}
\end{equation}

\subsection{Applying the rotations}
\label{CxAR}

The circuit of applying the rotations has expanding symmetry and is self-similar. This simplifies its implementation and complexity estimation. This also makes it a good candidate for pipelining. A high level look at the design of this circuit shows that it is made of 4 similar blocks each capable of applying single Given's Rotation. The critical path is determined by the circuit that applies double Given's rotations.

\begin{equation}
\begin{aligned}
\Delta_{Applying Rotations}= 2 \times \Delta_{One Given's Rotation} \\
= 2 \times ( \Delta_{\Lambda bit Multiply}+\Delta_{\Lambda bit Add}+ \Delta_{\Lambda bit Scale})
\end{aligned}
\label{eq:CxDAR}
\end{equation}

\begin{equation}
\begin{aligned}
A_{Applying Rotations}= 4 \times A_{One Given's Rotation} \\
= 4 \times ( \lambda \times A_{\Lambda bit Multiply}+ 4 \times A_{\Lambda bit Add}\\ + 4 \times A_{\Lambda bit Scale})
\end{aligned}
\label{eq:CxDAR}
\end{equation}

Assuming that the circuit of multiply as it represented in Figure~\ref{fig:Multiply} the followings equation will calculate the area and latency of the Multiply circuit.

\begin{equation}
\begin{aligned}
\Delta_{Multiply}= \Delta_{\lambda bit Add} + \Delta_{\Lambda bit B-Shifter} \\ + \Delta_{2 Inp XOR} + \Delta_{\Lambda bit Add} 
\end{aligned}
\label{eq:CxDMult}
\end{equation}

\begin{equation}
\begin{aligned}
A_{Multiply}= A_{\lambda bit Add}  + A_{\Lambda bit Add}  +\Lambda \times A_{2 Inp XOR} \\ +  A_{\Lambda bit B-Shifter} + (\Lambda + \lambda) \times A_{2 Inp Mux}
\end{aligned}
\label{eq:CxDMult}
\end{equation}

Considering the circuit represented in Figure~\ref{fig:Multiply} for Scaling, we calculate the delay and area requirement of this circuit for a 32 bit representation.\par 

\begin{equation}
\begin{aligned}
\Delta_{Scale}= \Delta_{32 bit B_Shifter} + \Delta_{32 bit Sub} \\ +  \Delta_{2 Inp Mux} + \Delta_{32 bit Add} +  \Delta_{2 Inp Mux} \\ + \Delta_{32 bit Add} +  \Delta_{2 Inp Mux}+ \Delta_{32 bit Add}\\ =\Delta_{32 bit B-Shifter} + 4 \times \Delta_{32 bit Add} + \Delta_{Inv} \\ + \Delta_{2 Inp Mux}+ 2 \time \Delta_{4 Inp Mux}\\ 
=\Delta_{32 bit B-Shifter} + 4 \times \Delta_{32 bit Add} \\ + \Delta_{Inv} + 5 \times \Delta_{2 Inp Mux}
\end{aligned}
\label{eq:CxDScale}
\end{equation}

\begin{equation}
\begin{aligned}
A_{Scale}= A_{32 bit B_Shifter} + A_{32 bit Sub} \\ + 32 \times A_{2 Inp Mux} + A_{32 bit Add} + 32 \times  A_{4 Inp Mux} \\ + A_{32 bit Add} + 32 \times A_{4 Inp Mux}+ A_{32 bit Add}\\ =A_{32 bit B-Shifter} + 4 \times A_{32 bit Add} + 32 \times A_{Inv} \\ + 32 \times A_{2 Inp Mux}+ 64 \times A_{4 Inp Mux}\\ 
=A_{32 bit B-Shifter} + 4 \times A_{32 bit Add} \\ +32 \times A_{Inv} + 224 \times A_{2 Inp Mux}
\end{aligned}
\label{eq:CxDScale}
\end{equation}

%
%



\section{Results}
We implemented the architecture presented in this work. Our implementation uses a 16 bit fixed point representation at the input (this number changes in the internal levels). We did not consider use of pipelining technique since pipelining and its achievable gain is orthogonal to the main idea of this paper.  The use of the pipelining and parallel hardware can be employed without any complications since the matrices are independent of each other. The use of parallel hardware can increase the throughput while it dose not effectively benefit the hardware efficiency and our design hardware efficiency will be poorer than some of the previous designs; on the other hand using pipelining will increase the throughput as well as hardware efficiency.\par

 The comparison of the presented method in this  with some of the state of the art works is presented in inverse chronological order in Table~\ref{tab:Comp}. This results show considerable improvement in energy per matrix target function over other designs. The sate of the art designs apply their algorithm on complex matrices, while our approach is employed to take advantage of reduced complexity achievable. This approach is valid since the recent publications are using SVD on the complex valued channel matrix of telecommunication systems; however for an application that requires the decomposition of a real valued matrix it is not that efficient. Any complex valued channel matrix can be converted to real valued channel matrices with four times the number of elements. To keep comparability we presented our synthesis results in a comparable form. A $2\times2$ matrix size in Table~\ref{tab:Comp} is a $2\times2$ complex valued matrix and equivalent to a $4\times4$ real valued matrix.\par

For comparison we considered the most resent designs that are able to calculate the SVD of nonsymmetric matrices and decompose a matrix to three matrices of $U,\ \Sigma,\ and\ V.$ Our goal is to show how our proposed design is able to provide energy efficient design with minimal hardware complexity. Due to the lack of pipelining and parallelism our proposed hardware does not achieve high hardware efficiency but the achievable gain from those methods is orthogonal to the benefit of this work original idea. In Table~\ref{tab:Comp} the results for various target functions are presented. In the telecommunication era power consumed to accomplish a task (as an indicator to show how fast the portable devices will drain their power using that particular device), and the throughput of a system are of the most importance. The value of the target function energy per matrix holds the effect of both important parameters. This parameter does not include the effect of hardware complexity. energy per matrix is the target function that can be used for comparison in this case. This function is able to project the effect of power consumption and throughput at the same time and ignores the effect of orthogonal techniques used for reducing the power, however, this function is still not able to eliminate the effect of hardware reuse in pipelining. \par

Beside the works compared in the Table~\ref{tab:Comp} authors in \cite{Comp4} use a method called supper linear SVD (SL-SVD) and are able to decompose a matrix sizes of $1 \times 1 \sim 4 \times 4$ very efficiently. The only downside to this algorithm is that it relies on the matrix quality and the channel characteristic. This matrix characteristic is harder to achieve with bigger size matrices. The convergence of this algorithm riles on all the singular values being different and in the case of two or more singular values being equal this algorithm never converges.\par  

In comparison with \cite{Comp1} our design achieves a lower throughput in smaller matrix sizes while our design without using any pipelining is able to achieve $23\%$ higher throughput for an $8\time8$ channel matrix. Authors in \cite{Comp1} use different bit sizes ($12 \sim 16$) for different matrix size ($2 \times 2 \sim 8 \times 8$) and uses different number of sweeps ($3 \sim 15$) for various matrix sizes. This method is also only designed to calculate the SVD of square matrices with even number of rows with complex elements. Our design on the other hand keeps the 16 bit\footnote{To achieve the required accuracy in 16 bit implementation in scaling circuit it is only critical to implement the B-Shifter, Subtractor, and first adder.} input accuracy for all the matrix sizes, and maps the number of sweeps used in \cite{Comp1} to the equivalent number for fast rotations method based on the values demonstrated in Figure~\ref{fig:NormDirectRelaxedSimp}. The proposed design is also able to decompose any square matrix of size $1 \times 1 \sim 8 \times 8$ (complex valued elements) and is able to as effectively decompose matrices with real valued element. In term of energy efficiency or our target function of energy per matrix our design provides $2.83 \sim 5.32$ (2.83 achieves from comparison of $8 \times 8$ matrices and 5.32 from comparison of $2 \times 2$ matrices) times better efficiency.\par

The authors of \cite{Comp6} and \cite{Comp5} do not provide the power consumption of their design.  The hardware complexity of the design represented in these works is considerably lower specially for lower sized matrices, however the throughput and normalized throughput of those designs is also lower. This in other word means that the hardware efficiency of these designs is lower and they are not good candidates for high throughput applications.\par

The authors in \cite{CompNew} use a Givens Rotation based design with a bipartite decomposition algorithm. First they convert a general matrix to a bidiagonalized matrix. The next step is to nullify the off diagonal elements. This design is capable of calculating both SVD and QRD (QR Decomposition). The proposed architecture in \cite{CompNew} is only capable of decomposing $4 \times 4$ matrices. Various techniques including pipelining, hardware sharing, and early termination are utilized to increase the hardware efficiency and throughput of the design. In the Table~\ref{tab:Comp} the power consumption of the design is mentioned with and without utilization of early termination process. this is to provide a fair comparison number since the gain achieved from early termination is application specific and also other designs could employ it. This design adopts a 12 bit implementation in contrast with our 16 bit implementation. The energy per matrix function of this design is $15.4\%$ better than our proposed design and its hardware complexity is $0.72\%$ higher.\par 

In comparison with \cite{Comp3} the energy per matrix of our proposed design is $649\sim36.131$ times for matrix sizes of $2\times 2 \sim 8 \times 8$ respectively. This benefit is achieved as the result of lower power usage of our work which is due to the simplicity of the hardware and eliminating the need for multiple pipeline registers. The achievable throughput and normalized achievable throughput are also higher in our proposed design.\par

In conclusion in term of Throughput as the target function, this design shows a superior performance compared to the works presented in \cite{Comp6},\cite{Comp3}, and \cite{Comp5}. The work in \cite{Comp1} is $4.51$ times better than our design for matrix size of $2\times2$; while this gap in the throughput result is reduced with the increase in the matrix size, and our design provides $1.24$ times better throughput for $8\times8$ matrices. This achievement is considerable since our design does not use any parallelism or pipelining. In term of hardware efficiency or Normalized Throughput function the design presented in this work provides better results than the work presented in \cite{Comp6},\cite{Comp3}, and \cite{Comp5}. While the normalized throughput of this work is $1.47\sim3.931$ times for matrix sizes of $2\times2\sim8\times8$ respectively compared to the work presented in \cite{Comp1}, this is not far from expected since we expected the hardware efficiency of our design to be poorer.

\begin{table*}[]
\centering
\caption{Comparison of different decomposition algorithms.}
\label{tab:Comp}
\tiny
\begin{threeparttable}
\begin{tabular}{ccccccc}
\hline
                                                           & This Work                                                                                                                                & \cite{Comp1}                                                   & \cite{Comp6}                    & \cite{CompNew}                               & \cite{Comp3}                                                   & \cite{Comp5}                                   \\ \hline \hline
Technology (nm)                                     & 90                                      & 90                                                               & 90                                                               & 90                  & 90                                                               & 90                                               \\ \hline
Algorithm                                           & ERFHSVD                                 & napSVD                                                           & 2-Sided Jacobi                                                   & GR Based            & Adaptive SVD                                                     & GK                                               \\ \hline
Functional                                          & U, $\Sigma$, and V                      & U, $\Sigma$, and V                                               & U, $\Sigma$, and V                                               & U, $\Sigma$, and V  & U, $\Sigma$, and V                                               & U, $\Sigma$, and V$\slash$Q, R, and P            \\ \hline
Matrix Size                                         & 1$\times$1 $\sim$ 8$\times$8\tnote{1} & 2$\times$2$\slash$4$\times$4$\slash$6$\times$6$\slash$8$\times$8 & 2$\times$2$\slash$4$\times$4$\slash$6$\times$6$\slash$8$\times$8 & 4$\times$4          & 1$\times$1$\slash$2$\times$2$\slash$3$\times$3$\slash$4$\times$4 & 4$\times$4$\slash$8$\times$8$\slash$16$\times$16 \\ \hline
Max. Frequency (MHz)                                & 125                                     & 752                                                              & 112                                                              & 143                 & 101.2                                                            & 400                                              \\ \hline
Gate Count (KGE)                                    & $116.7$\slash$448.3$\slash\$1756.1      & 359                                                              & 378                                                              & 451.53              & 543.9                                                            & 54.5                                             \\ \hline
Power (mW)                                          & $16.8 $\slash$64.5$\slash\$252.8        & 402$\slash$595$\slash$673$\slash$770                             & -                                                                & 164.4$\slash$218.29\tnote{2} & 125                                                              & -                                                \\ \hline
Throughput (MMatrices$\slash$s)                     & $41.7 $\slash$8.93$\slash\$2.08         & 188.1$\slash$15.7$\slash$6.3$\slash$1.68                         & -$\slash$0.229$\slash$-$\slash$0.09866                           &             35.75        & 0.4791                                                           & 0.01391$\slash$0.00136$\slash$-                  \\ \hline
Normalized Throughput (Matrices$\slash$s$\slash$GE) & $357$\slash$19.92$\slash\$1.19          & 524$\slash$43.7$\slash$17.5$\slash$4.68                          & -$\slash$0.607$\slash$-$\slash$0.261                             &      7.87                & 0.88                                                             & 0.255$\slash$0.025$\slash$-                      \\ \hline
Energy per Matrix (nJ)                              & 0.4021$\slash$7.222$\slash$121.4      & 2.14$\slash$38$\slash$107.4$\slash$343.6                         & -                                                                & 4.6$\slash$6.11\tnote{2}                     & 260.9                                                            & -                                                \\ \hline
\end{tabular}
\begin{tablenotes}
\item[1]{Only the results for 2$\times$2$\slash$4$\times$4$\slash$8$\times$8 is mentioned in the table in favor of simplifying the presentation.}
\item[2]{The power numbers are considered with$\slash$without early termination method.}
\end{tablenotes}
\end{threeparttable}
\end{table*}

\section{Conclusion}
\label{Concl}

In this work, for the first time we presented an algorithm that directly estimates the fast rotations for singular value decomposition of a non symmetric matrix. 
This method, unlike the previous efforts to implement an eigenvalue decomposition, is able to provide the "Normalized" results.
An implementation is presented for $2 \times 2$ matrix as the basic block cell of any matrix of higher size. 
Unlike the previous efforts the proposed hardware does not require any floating point representation or hardware.
The analysis of implementation requirement and complexity is also presented. 
The hardware complexity of the $2\times2$ matrix decomposer is much lower than the previous floating-point implementations.
This design provides $2.83\sim649$ times better energy per matrix performance compared to the most resent designs.


%

%
%
%
%
%

\ifCLASSOPTIONcaptionsoff
  \newpage
\fi



%

\bibliography{MyBib2}
\bibliographystyle{IEEETran}

%

%
%
%




\end{document}